\documentclass{alex-preprint}
\usepackage{amsthm}

\newcommand{\ShockMin}{u_\ell}
\newcommand{\ShockMax}{u_r}
\newcommand{\ds}{\displaystyle}

\definecolor{darkmagenta}{rgb}{0.55, 0.0, 0.55}
\newtheorem{theorem}{Theorem}

\newtheorem{prop}{Proposition}
\crefname{prop}{Proposition}{Propositions} 

\title{Shock selection in reaction--diffusion equations with partially negative diffusivity using nonlinear regularisation}
\author[1]{Thomas Miller}
\author[1]{Alexander K. Y. Tam}
\author[2]{Robert Marangell}
\author[2]{Martin Wechselberger}
\author[1,*]{Bronwyn H. Bradshaw-Hajek}
\affil[1]{UniSA STEM, The University of South Australia, Mawson Lakes SA 5095, Australia}
\affil[2]{School of Mathematics and Statistics, The University of Sydney, Sydney NSW 2006, Australia}
\corrauthor[*]{\href{mailto:bronwyn.hajek@unisa.edu.au}{bronwyn.hajek@unisa.edu.au}}

\abstract{We consider a general reaction--nonlinear-diffusion equation with a region of negative diffusivity, and show how a nonlinear regularisation selects a shock position. Negative diffusivity can model population aggregation, but leads to shock-fronted solutions for population density. In general the shock position is non-unique. Previous studies have defined shock selection criteria such as the equal area rule, and shown how these arise from specific regularisations to the reaction--diffusion equation. In this work, we show that a nonlinear regularisation leads to travelling wave solutions where the shock is selected according to a modified equal area rule. Adjusting the nonlinearity in the regularisation moves the shock location. We focus on attaining shocks that conserve diffusivity across the shock, and demonstrate that this condition yields the longest possible shock length. Using geometric singular perturbation theory, we prove the existence of shock-fronted travelling wave solutions with continuous diffusivity, show how to construct them, and demonstrate that they correspond to a unique wave speed. Numerical solutions align with theoretical predictions for shock position and wave speed, confirming that a single regularisation term can vary the shock position and attain shocks with continuous diffusivity.}


\begin{document}
    \onehalfspacing
    \maketitle
    \section{Introduction}\label{section:Introduction}
        In this paper, we establish a shock selection rule based on continuous diffusivity for a reaction--diffusion equation with a region of negative diffusivity using a nonlinear regularisation. Reaction--diffusion equations are partial differential equations that are widely used in mathematical biology and population dynamics~\parencite{murray2002mathematical,edwards2018compactly,el2021invading,gatenby1996reaction,lewis1993allee}. In one spatial dimension the general reaction--diffusion equation is
        \begin{equation}
            \label{equation:reaction_diffusion_equation}%
            \pd{u}{t} = \pd{}{x}\left(D(u) \pd{u}{x}\right) + R(u) = \pdn{\Phi(u)}{x}{2} + R(u),
        \end{equation}
        where \(u(x,t)\) is the population density or concentration as a function of space \(x\) and time \(t\), \(D(u)\) is the diffusivity, and \(R(u)\) is the reaction term describing local production or loss of \(u.\) In this work, it will often be useful to write~\Cref{equation:reaction_diffusion_equation} in terms of the flux potential (or potential), \(\Phi(u) = \int D(u)\df{u}.\)
    
        The diffusivity \(D(u)\) describes how the population moves. Diffusion according to Fick's second law obeys \(D(u) = D,\) where \(D\) is a positive constant. This scenario is linear diffusion. However, in many applications the diffusivity \(D(u)\) is non-constant and varies with \(u\). This situation is nonlinear diffusion. Nonlinear diffusivity usually satisfies \(D(u) \geq 0\) for all \(u,\) which leads to the population spreading from areas of high density to areas of low density~\parencite{murray2002mathematical}. However, in certain contexts individuals in the population might prefer to move from areas of low density to areas of high density. For example, individuals might aggregate to avoid extinction~\parencite{grindrod1988models}, and animals aggregating for social reasons is common in ecology~\parencite{grunbaum1994modelling}. Reaction--diffusion equations can model aggregation by allowing \(D(u) < 0\) for some values of \(u.\)  
            
        Negative diffusivity models arise in the continuum limit of discrete models that explicitly include aggregation. For example, \textcite{turchin1989population} derived a diffusion equation \(\left(R(u) = 0\right)\) with a region of negative diffusivity by considering a random walk where individuals at low density were more likely to be attracted to each other, whereas at high density the attraction was reduced or even turned to repulsion. \textcite{johnston2017co} derived reaction--diffusion equations with positive-negative-positive quadratic diffusivities from an agent based model where isolated agents and group agents undergo different rates of birth, death and movement. The resulting reaction--diffusion equations have been studied analytically. Like well-known examples with \(D(u) \geq 0\)~\parencite{murray2002mathematical}, the existence of travelling wave solutions to models with negative diffusivity has been well-established~\parencite{ferracuti2009travelling,kuzmin2011front,li2020travelling}.

\begin{figure}[htbp!]
            \centering
            \subcaptionbox{Diffusivity, \(D(u).\)\label{figure:diffusion_sketch}}{\includegraphics[width=0.49\linewidth]{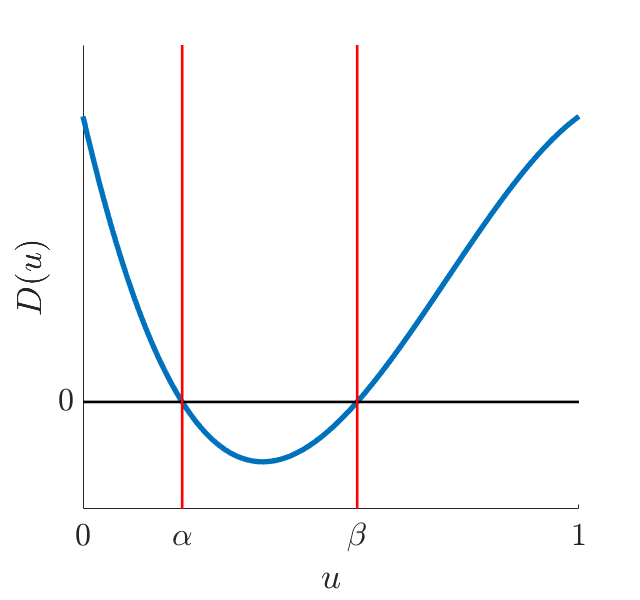}}
            \subcaptionbox{Potential, \(\Phi(u).\)\label{figure:phi_shocks}}{\includegraphics[width=0.49\linewidth]{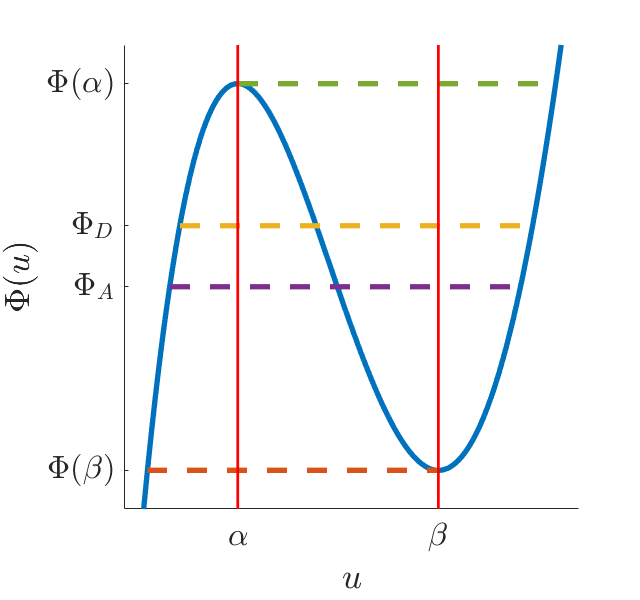}}
            \caption{A decreasing-increasing diffusivity, \(D(u)\), and corresponding potential \(\Phi(u)\). (a) Diffusivity \(D(u),\) such that \(D(u)>0\) for \(u \in [0,\alpha)\bigcup(\beta,1]\), \(D(u)<0\) for \(u \in (\alpha,\beta)\), and \(D(u)=0\) when \(u = \alpha, u = \beta\). (b) The potential, \(\Phi(u)\), showing four possible shock positions (dashed horizontal lines - colour online). Green: Shock from the upper knee of \(\Phi,\) where \(\ShockMin = \alpha.\) Yellow: Shock satisfying continuous diffusivity. Purple: Shock satisfying the equal area in \(\Phi(u)\) rule. Red: Shock at the lower knee of \(\Phi,\) where \(\ShockMax = \beta\). Vertical lines (red online) highlight the zeros of the diffusivity. Between the solid vertical lines, the diffusivity is negative and \(\Phi(u)\) is decreasing.}
            \label{figures:diffusion_sketch+phi_shocks}
        \end{figure}

        Reaction--diffusion equations with negative diffusivity give rise to shock solutions~\parencite{witelski1995shocks,pego1989front,pego2007lectures,pesce2021degenerate}, whereby the solution is sharp-fronted or, in the limiting case, discontinuous. We consider diffusivities that are only negative for a region \(u \in (\alpha,\beta)\) where \(u = \alpha,\) \(u = \beta\) are the zeros of the diffusivity and \(0 < \alpha < \beta < 1\), see for example \Cref{figure:diffusion_sketch}. With \(D(u) < 0\) in a region \(u \in (\alpha,\beta)\), shock solutions become possible where a solution jumps from one density \(\ShockMax \geq \beta\) instantaneously to another density \(\ShockMin \leq \alpha\) (or vice versa), such that \(u\) never takes the values in the region \((\ShockMin,\ShockMax)\) \(\supset (\alpha, \beta)\) and hence completely avoids the region of negative diffusivity altogether~\parencite{witelski1995shocks}.
        
        The location of a shock is defined by the values of the density at its endpoints, \(u = \ShockMin \leq \alpha\) and \(u = \ShockMax \geq \beta\). However, in terms of existence, these endpoints, and hence the shock location, are in general not unique. \Cref{figure:phi_shocks} shows the potential \(\Phi(u)\) corresponding to the diffusivity in~\Cref{figure:diffusion_sketch}, with four possible shock positions illustrated with dashed lines. The diffusivity and potential are related by \(D(u) = \Phi'(u)\) so that when \(\Phi(u)\) is increasing the diffusivity is positive, while when \(\Phi(u)\) is decreasing the diffusivity is negative. The local maximum and minimum of \(\Phi(u)\) occur where the diffusivity is zero. Shock solutions to reaction--diffusion equations avoid regions of negative diffusivity, connecting the regions of increasing \(\Phi(u)\) and skipping over the regions of decreasing \(\Phi(u)\). 
        
        To specify the values of the shock endpoints uniquely, and hence determine the shock solution profile, two conditions are required to supplement the reaction--diffusion equation \eqref{equation:reaction_diffusion_equation}. The first condition is that \(\Phi(u)\) is continuous across the shock, as per all shocks illustrated in~\Cref{figure:phi_shocks}. Continuous \(\Phi(u)\) means that
        \begin{equation}
            \label{equation:continous_phi_shock_condition}%
            \Phi(\ShockMin) = \Phi(\ShockMax) = \Phi_S,        
        \end{equation}
        where \(\Phi_S\) is the value of \(\Phi(u)\) that is conserved across the shock. Since \(\Phi(u)\) is continuous across the shock, \(\Phi(\beta) \leq \Phi_S \leq \Phi(\alpha)\). The continuous \(\Phi\) condition is necessary (but not sufficient) 
        for the solution with the shock inserted to be a weak solution to the governing reaction--diffusion equation~\parencite{whitham2011linear,witelski1995shocks}. Unless explicitly stated otherwise, in this paper we will always assume that $\Phi$ is continuous accross the shock. 
        
\begin{figure}[htbp!]
            \centering
            \subcaptionbox{\label{figure:equal_area_in_phi_shock_condition}}{\includegraphics[width=0.49\linewidth]{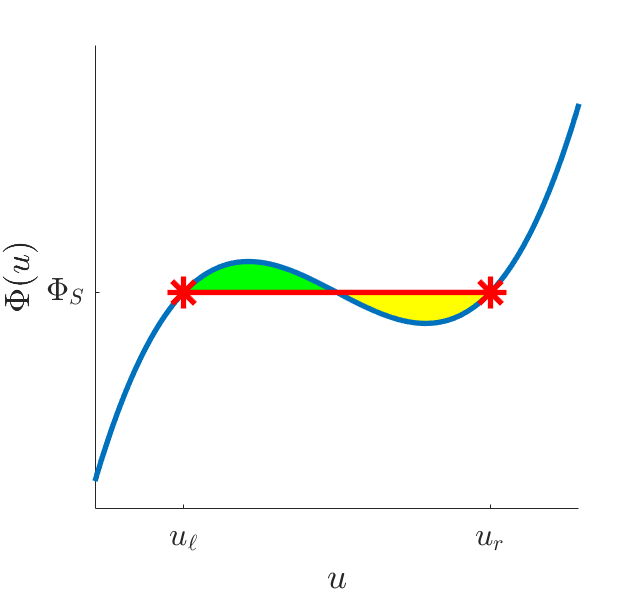}}
            \subcaptionbox{\label{figure:continous_diffusion_shock_condition}}{\includegraphics[width=0.49\linewidth]{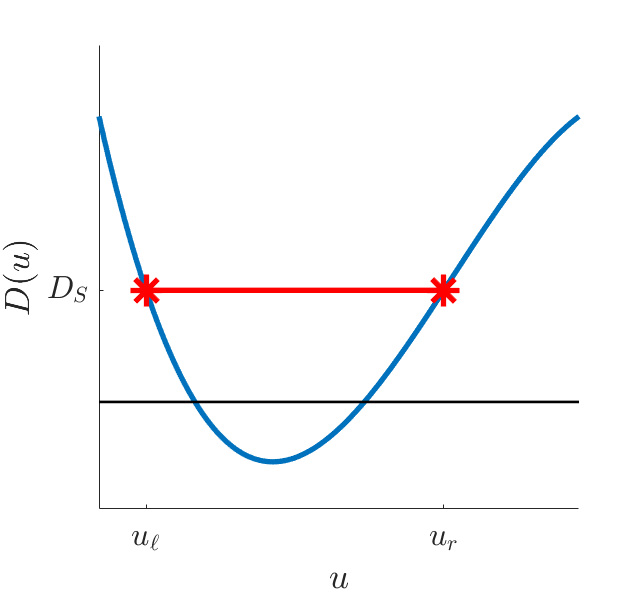}}
            \caption{Two possible shock selection rules. (a) Equal area in \(\Phi(u)\) rule. Red stars indicate the shock endpoints \(\ShockMin\) and \(\ShockMax,\) and the red horizontal line indicates the value of \(\Phi_S\) that is conserved across the shock. The areas of the two shaded regions are equal. (b) Continuous \(D(u)\) rule. Red stars indicate the shock endpoints, \(\ShockMin\) and \(\ShockMax\), and the red horizontal line shows that the diffusivity is the same at both \(\ShockMin\) and \(\ShockMax\).}
        \end{figure}

        A second commonly used condition is the equal area in \(\Phi\) rule~\parencite{pego1989front,goh2018dynamics,popescu2004model}
        \begin{equation}
            \label{equation:equal_area_in_phi_shock_condition}%
            \int_{\ShockMin}^{\ShockMax} \left(\Phi(u) - \Phi_S\right) \df{u} = 0,  
        \end{equation}
        which is illustrated in~\Cref{figure:equal_area_in_phi_shock_condition} and is represented by the dotted line, second from the bottom (purple online) in \Cref{figure:phi_shocks}. When the diffusivity is symmetric about its vertex, the shock position given by the equal area in \(\Phi(u)\) rule also satisfies continuous diffusivity across the shock, \emph{i.e.}
        \begin{equation}
            \label{equation:continous_diffusion_shock_condition}%
            D(\ShockMin) = D(\ShockMax) = D_S,
        \end{equation}
        which is illustrated in~\Cref{figure:continous_diffusion_shock_condition}. An example of a symmetric diffusivity is a quadratic. However, when the diffusivity is non-symmetric about its vertex, application of the equal area rule does not result in the diffusivity being continuous across the shock (see \cite{miller2024analytic}). \Cref{figure:phi_shocks} shows a shock location with a continuous diffusivity as a yellow dotted line in the non-symmetric case. 

        Shock selection conditions can correspond to regularising the reaction--diffusion equation~\parencite{witelski1995shocks}. Regularisation involves adding higher-order terms of small size (\(\mathcal{O}(\varepsilon)\) or smaller, as \(\varepsilon \to 0\)) to the right-hand side of the reaction--diffusion equation~\cref{equation:reaction_diffusion_equation}. Specific regularisations correspond to specific shock solutions from among those available. For diffusion and reaction--diffusion equations the equal area in \(\Phi(u)\) rule is associated with a fourth order non-local regularisation term. \textcite{witelski1995shocks} demonstrated this for diffusion equations by analysing the \(\varepsilon \rightarrow 0\) limit of the Cahn-Hilliard equation,
        \begin{equation}
            \label{equation:cahn-hilliard}%
            \pd{u}{t} = \pdn{}{x}{2}\left(\Phi(u) - \varepsilon^2 \pdn{u}{x}{2}\right),
        \end{equation}
        using perturbation methods. \textcite{li2021shock} proved the existence of shock solutions to reaction--diffusion equations by applying geometric singular perturbation theory (GSPT) to a fourth order regularised reaction--diffusion equation
        \begin{equation}
            \label{equation:standard_fourth_order_regularisation}%
            \pd{u}{t} = \pdn{}{x}{2}\left(\Phi(u) - \varepsilon^2 \pdn{u}{x}{2}\right) + R(u),
        \end{equation}
        and a third order viscous regularised reaction--diffusion equation
        \begin{equation}
            \label{equation:regularised_reaction_diffusion_mixed}%
            \pd{u}{t} = \pdn{}{x}{2}\left(\Phi(u) + \varepsilon \pd{u}{t}\right) + R(u),  
        \end{equation}
        in the \(\varepsilon \rightarrow 0\) limit. The third-order viscous regularisation~\cref{equation:regularised_reaction_diffusion_mixed} leads to shock positions that either jump to or from one of the zeros of the diffusivity~\parencite{li2021shock,witelski1996structure} (as per the red or green dotted lines in~\Cref{figure:phi_shocks}). By combining both the non-local and viscous regularisation terms
        \begin{equation}
            \label{equation:mixed_regularisation_and_fourth_order_regularisation}%
            \pd{u}{t} = \pdn{}{x}{2}\left(\Phi(u) + \omega \varepsilon \pd{u}{t} - \varepsilon^2 \pdn{u}{x}{2}\right) + R(u),
        \end{equation}
        the position of the shock can be controlled by using \(\omega\) to balance the relative weighting of the two regularisation terms~\parencite{bradshaw2023geometric}. A similar result was found for a convective viscous Cahn-Hilliard equation~\parencite{witelski1996structure}.
    
       The aim of this paper is to show that a nonlinear regularisation term in a reaction--diffusion equation \eqref{equation:reaction_diffusion_equation} can be used to specify a shock position the preserves the diffusivity \(D(u)\) across the shock. In \Cref{section:Statement_of_Results} we state our main results, which we prove in~\Cref{section:Results_and_Discussion}. We derive using a perturbation method and by using geometric singular perturbation theory a modified equal area rule~\Cref{equation:modified_equal_area_rule} which arises from a nonlinear fourth order regularisation~\Cref{equation:nonlinear_fourth_order_regularisation_2_2} which will allow us to control the position of the shock. We explore shock solutions numerically in~\Cref{section:numerical}, to show how using a nonlinear fourth order regularisation enables control of the position of a shock.

    \section{Statement of Results}\label{section:Statement_of_Results}

        Unless otherwise stated, we assume that the diffusivity, \(D(u)\) has two zeros \(\alpha < \beta\) between \(u=0\) and \(u=1\), (\(D(\alpha)=D(\beta)=0\)). The diffusivity is negative, \(D(u) < 0\), for \(\alpha < u < \beta,\) and positive for \(0 < u < \alpha\) and \(\beta < u < 1.\) If the diffusivity decreases monotonically for \(0\le u\le\alpha\) and increases monotonically for \(\beta\le u\le1\) we refer to it as decreasing-increasing. In addition, {\it symmetric diffusivity} refers to a diffusivity that is symmetric about the midpoint of its zeros. A non-symmetric diffusivity is not symmetric about the midpoint of its zeros. As shown by~\textcite{miller2024analytic} if the diffusivity \(D(u)\) is symmetric about the midpoint of its zeros the shock position satisfying the equal area in \(\Phi(u)\) rule also satisfies a continuous diffusivity rule, while for non-symmetric diffusivity this is not the case in general. We illustrate this further in~\Cref{section:nonsymmetric_diffusion}, before proceeding to prove the following results. 
            
        \begin{prop} \label{lemma:existence_of_shock_position_equal_area}
            For a continuous \(D(u)\) where \(D(u) < 0\) for \(u \in (\alpha,\beta)\) and \(D(u) \geq 0\) for \(u \notin (\alpha,\beta)\), it is always possible to find a shock position (\(\ShockMax\), \(\ShockMin\)) such that $\Phi(u)$ is continuous across the shock and that satisfies the equal area in \(\Phi(u)\) rule. That is:
            \begin{enumerate}
                \item \(\Phi(\ShockMax) = \Phi(\ShockMin) = \Phi_S\), and 
                \item  \(\ds\int_{\ShockMin}^{\ShockMax} \Phi(u) - \Phi_S\, \df{u} = 0.\)
            \end{enumerate}                        
        \end{prop}

        \begin{prop} \label{lemma:existence_of_shock_position_equal_diffusion}
            For a continuous \(D(u)\) where \(D(u) < 0\) for \(u \in (\alpha,\beta)\) and \(D(u) > 0\) for \(u \notin (\alpha,\beta)\), it is always possible to find a shock position (\(\ShockMax\), \(\ShockMin\)) such that $\Phi(u)$ is continuous across the shock and that satisfies a continuous diffusivity rule. That is:
            \begin{enumerate}
                \item \(D(\ShockMax) = D(\ShockMin) = D_S\), and 
                \item  \(\Phi(\ShockMax) = \Phi(\ShockMin) = \Phi_S\).
            \end{enumerate}
            Further, when \(D(u)\) is decreasing-increasing, the shock position satisfying continuous diffusivity is unique.            
        \end{prop}        
In light of \cref{lemma:existence_of_shock_position_equal_diffusion}, we will call a shock satisfying the conditions of \cref{lemma:existence_of_shock_position_equal_diffusion}, a {\em continuous diffusivity shock}.
        
        Proofs of~\Cref{lemma:existence_of_shock_position_equal_area,lemma:existence_of_shock_position_equal_diffusion} are given in~\Cref{proof:existence_of_shock_positions}. When the diffusivity is decreasing-increasing, the shock determined by requiring that the diffusivity is continuous across the shock is unique and, as we shall see (\Cref{lemma:longest_shock}), coincides with the position that gives the longest possible jump in density.

        \begin{prop}\label{lemma:longest_shock}
            For a continuous \(D(u)\) which is decreasing-increasing, the shock position satisfying continuous diffusivity is the longest possible shock. That is, the continuous \(D(u)\) shock maximises shock length,
            \begin{equation}
                \label{equation:shock_length}%
                S_L = \ShockMax - \ShockMin.
            \end{equation}
        \end{prop}

        This is proved in \Cref{proof:longest_shock} and is illustrated in \Cref{figure:shock_lengths_delta_0.5}. Our main result is the following, concerning the location and uniqueness of the shock.

        \begin{theorem}\label{lemma:perturbation}             
            When travelling wave solutions with shocks exist, shock positions satisfy a modified equal area rule,
            \begin{equation}
                \label{equation:modified_equal_area_rule}%
                \int_{\ShockMin}^{\ShockMax} \frac{\Phi(u) - \Phi_S}{f(u)} \df{u} = 0,
            \end{equation}
            and correspond to the \(\varepsilon \rightarrow 0\) limit of a nonlinear fourth order regularised reaction--diffusion equation,        
            \begin{equation}                    
                \label{equation:nonlinear_fourth_order_regularisation_2_2}%
                \pd{u}{t} = \pdn{}{x}{2}\left(\Phi(u) - \varepsilon^2 f(u) \pdn{u}{x}{2}\right) + R(u),
            \end{equation}
            where the diffusivity \(D(u)\) is negative for \(u \in (\alpha,\beta)\), and positive for \(u \notin (\alpha,\beta)\) and \(f(u) > 0\). This holds for arbitrary \(R(u)\), including \(R(u) = 0\).
        \end{theorem}
     
      \Cref{lemma:perturbation} says that a modified equal area rule that directly corresponds to a choice of regularisation can be used to specify a shock position. In particular we will use it to specify a shock position that satisfies the continuous diffusivity rule \eqref{equation:continous_diffusion_shock_condition}. \Cref{section:perturbation} contains the proof. In terms of existence of shock-fronted solutions to \cref{equation:reaction_diffusion_equation} with cubic reaction terms, we rely on standard results from geometric singular perturbation theory.

        \begin{theorem}\label{theorem:GSPT}
            For the cubic reaction \(R(u) = u(1 - u)\left(u - \gamma\right)\) with \(\gamma \in (0,1)\), geometric singular perturbation theory (GSPT) shows that travelling wave solutions with a shock satisfying the modified equal area rule~\Cref{equation:modified_equal_area_rule} exist in the \(\varepsilon \rightarrow 0\) limit of a nonlinear fourth order regularised reaction--diffusion equation,~\Cref{equation:nonlinear_fourth_order_regularisation_2_2}, where the diffusivity \(D(u)\) is negative for \(u \in (\alpha,\beta)\), and positive for \(u \notin (\alpha,\beta)\) with \(0 < \alpha < \beta < 1\), and \(f(u) > 0\).
        \end{theorem}

        This is proved in~\Cref{section:GSPT} by writing the regularised reaction--diffusion equation~\Cref{equation:nonlinear_fourth_order_regularisation_2_2} as a fast-slow system and applying well-known techniques of the theory. Our next result uses~\Cref{theorem:GSPT,lemma:perturbation} and exponential \(f(u)\) to prove that we can use the nonlinear regularisation~\cref{equation:nonlinear_fourth_order_regularisation_2_2} and the modified equal area rule to obtain a travelling wave solutions with continuous diffusivity shock.

        \begin{theorem}\label{lemma:finding_A}    
            For a polynomial diffusivity \(D(u)\) with \(D(u) < 0\) for \(u \in (\alpha,\beta)\), \(D(u) > 0\) for \(u \notin (\alpha,\beta)\), a given shock position \(\ShockMax, \ShockMin\), and with \(f(u) = \e^{-Au}\) as in \cref{equation:nonlinear_fourth_order_regularisation_2_2}, it is always possible to find a value of \(A\) such that the modified equal area rule~\Cref{equation:modified_equal_area_rule} is satisfied and yields the given shock position when \(D(\ShockMax) > 0\) and \(D(\ShockMin) > 0\).
        \end{theorem}

        \Cref{section:Exponential_Regularisation} contains the proof of~\Cref{lemma:finding_A}. In~\Cref{figure:A_values} we show how the value of \(A\) corresponding to shocks satisfying the continuous diffusivity rule \eqref{equation:continous_diffusion_shock_condition} varies with the degree of asymmetry in the diffusivity.

    \section{Proofs and Discussion}\label{section:Results_and_Discussion}
        \subsection{Equal area and continuous diffusivity shocks do not coincide for non-symmetric diffusivity}\label{section:nonsymmetric_diffusion}        
            The equal area in \(\Phi(u)\) and  continuous diffusivity shocks only coincide~\parencite{miller2024analytic} if the diffusivity \(D(u)\) is symmetric about the midpoint of its zeros. To demonstrate this we use the diffusivity
            \begin{equation}
                \label{equation:quadratic_cubic_diffusion}%
                D(u) = (u - a)(u - b  - \delta u^2),
            \end{equation}
            that depends on the parameters \(a,\) \(b,\) and \(\delta\). When \(\delta = 0\), the diffusivity is quadratic with zeros at \(u = a\) and \(u = b\), and is symmetric about its vertex which lies at the mid-point of its zeros. If \(\delta \neq 0\), the diffusivity~\Cref{equation:quadratic_cubic_diffusion} is cubic and not symmetric, and has zeros at \(u = a\) and \(u = (1 \pm \sqrt{1 - 4b\delta})/{2\delta}\). We choose \(a,\) \(b,\) and \(\delta\) such that there are exactly two zeros for \(u \in (0,1)\) and such that $D(0)>0$ and $D(1)>0$. We label the zeros \(\alpha\) and \(\beta\) such that \(0 < \alpha < \beta < 1\). The potential, \(\Phi(u)\), corresponding to \Cref{equation:quadratic_cubic_diffusion} is
            \begin{equation}
                \label{equation:phi_for_quadratic_cubic_diffusion}%
                \Phi(u) = \int_{u^*}^u D(\Bar{u}) \df{\Bar{u}} = \frac{1}{12}u\left(-3\delta u^3 + 4a\delta u^2 + 4u^2 - 6au - 6bu + 12ab\right) + c^*,
            \end{equation}
            where \(u^*\) is an arbitrary base point that determines \(c^*\), the constant of integration. Without loss of generality, we take \(u^*=c^* = 0\).
            \begin{figure}[htbp!]
                \centering
                \includegraphics[width=0.65\linewidth]{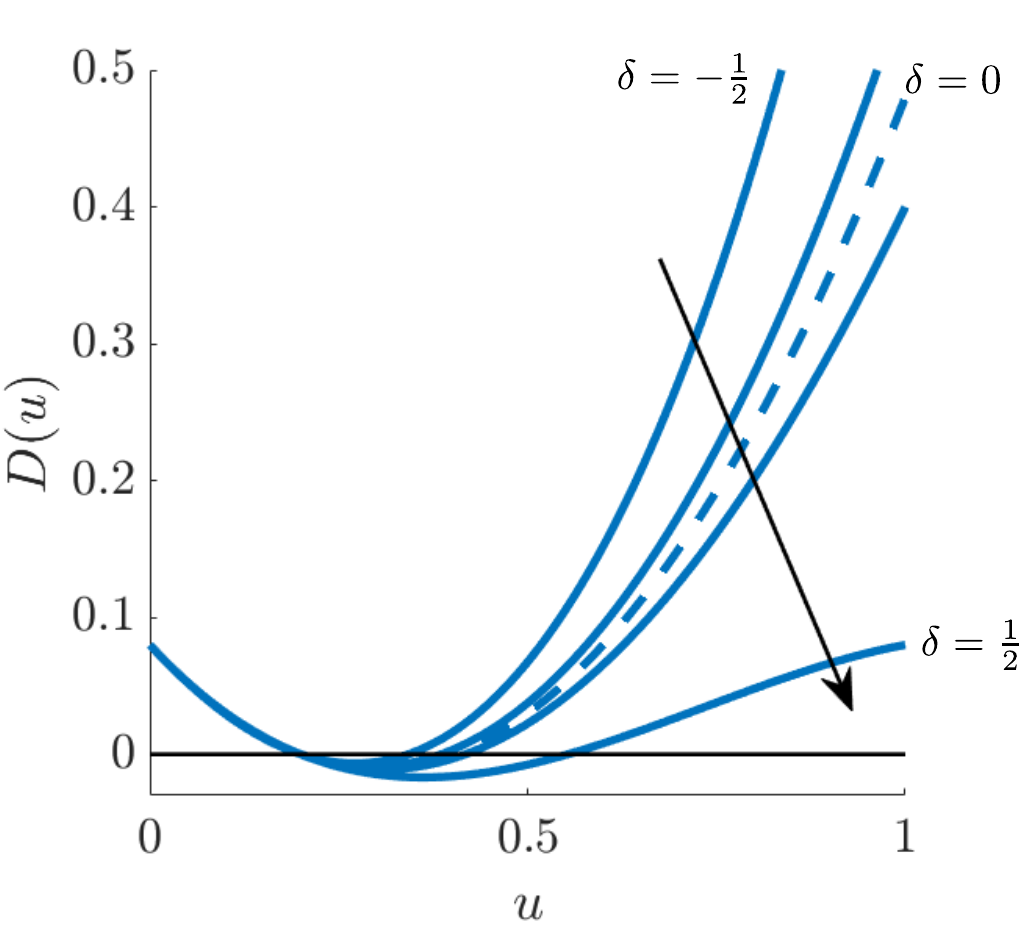}
                \caption{Quadratic/cubic diffusivities~\cref{equation:quadratic_cubic_diffusion}, with \(a = 0.2,\) \(b = 0.4\), and varying \(\delta.\) The dashed curve is quadratic \((\delta = 0)\). The other curves have \(\delta\neq0\) and are not symmetric. The diffusivities are plotted for \(\delta \in \{-0.5, -0.1, 0, 0.1, 0.5\},\) and the arrow indicates the direction of increasing \(\delta\).}
                \label{figure:Diffusion_multiple_delta}
            \end{figure}
    
\begin{figure}[htbp!]
                \centering
                \subcaptionbox{\label{figure:shock_min_values}}{\includegraphics[width=0.49\linewidth]{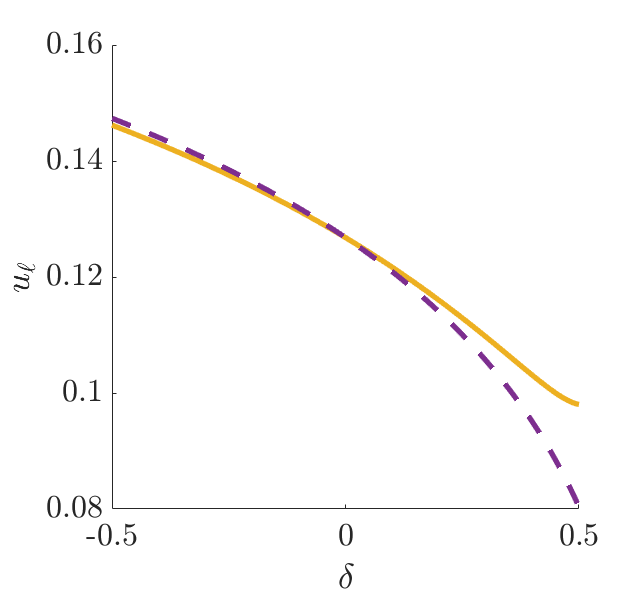}}
                \subcaptionbox{\label{figure:shock_max_values}}{\includegraphics[width=0.49\linewidth]{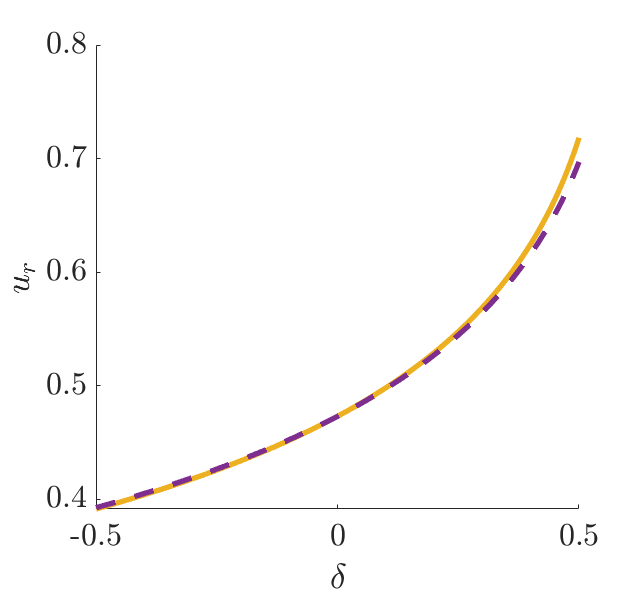}}
                \caption{Shock end-points for the quadratic/cubic diffusivity~\cref{equation:quadratic_cubic_diffusion}, with \(a = 0.2\), \(b = 0.4,\) and varying \(\delta.\) Equal area in \(\Phi(u)\) shock positions are shown in dashed purple, and continuous diffusivity shock position are the solid yellow curve. (a) Lower/left endpoint, \(\ShockMin\). (b) Upper/right endpoint, \(\ShockMax\).}
                \label{figures:shock_position_values}
            \end{figure}

            With \(\delta\neq 0\) the equal area in \(\Phi(u)\) shock position and the continuous diffusivity shock position do not in general coincide. The two different shock locations can ostensibly be found by solving equations \eqref{equation:continous_phi_shock_condition} and \eqref{equation:equal_area_in_phi_shock_condition} or equations \eqref{equation:continous_phi_shock_condition} and \eqref{equation:continous_diffusion_shock_condition} respectively. For even the lowest degrees this is analytically tedious, but \cref{lemma:existence_of_shock_position_equal_area,lemma:existence_of_shock_position_equal_diffusion} show that this can always be accomplished. \Cref{figures:shock_position_values} shows how the endpoints of the shock change with \(\delta\) for each shock position, and demonstrate that the shocks coincide when \(\delta=0\). Throughout this work, we consider diffusivities with \(\delta \neq 0\) to distinguish between equal area in \(\Phi\) and continuous diffusivity shocks.      
    
        \subsection{Existence of equal area and continuous diffusivity shock positions} \label{proof:existence_of_shock_positions}        
            In this subsection we prove~\Cref{lemma:existence_of_shock_position_equal_area,lemma:existence_of_shock_position_equal_diffusion}. These establish the existence of shocks satisfying the equal area in \(\Phi\) and continuous diffusivity rules.                  
                                
\begin{figure}
                \centering
                \includegraphics[width=0.65\linewidth]{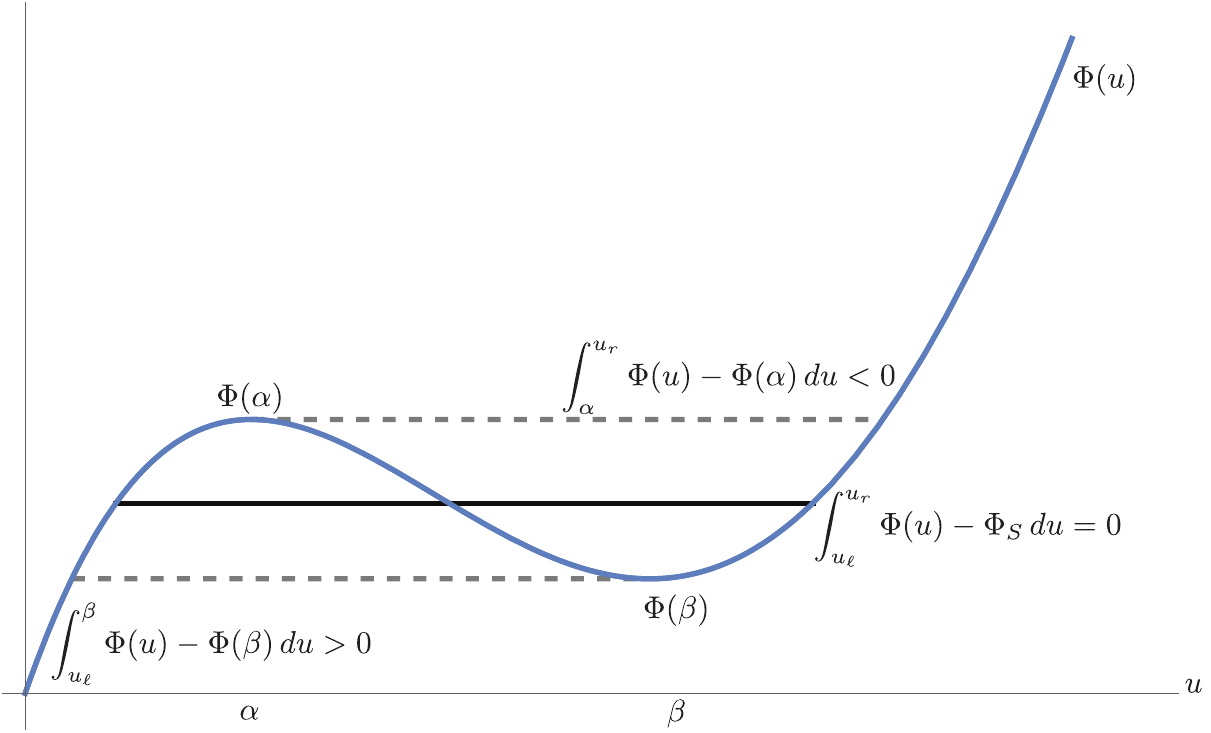}
                \caption{A schematic of the proof of \Cref{lemma:existence_of_shock_position_equal_area}.}
                \label{fig:lemma1}
            \end{figure}

            \begin{proof}[Proof of \Cref{lemma:existence_of_shock_position_equal_area}: existence of an equal area in $\Phi(u)$ shock]  
                The intermediate value theorem can be used to argue that it is always possible to specify a shock location that satisfies the equal area in \(\Phi(u)\) rule as follows. If the lower endpoint of the shock coincides with the first zero of the diffusivity, \(\ShockMin = \alpha\), then \(\Phi_S = \Phi(\alpha)\) is the local maximum of the potential. For all \(u \in (\ShockMin,\ShockMax)\), \(\Phi(u) < \Phi_S\), so that the left hand side of \Cref{equation:equal_area_in_phi_shock_condition}, the quantity $\int_\alpha^{u_r} \Phi(u) - \Phi_S \, du$ is negative. Conversely, when the upper endpoint of the shock coincides with the second zero of the diffusivity, \(\ShockMax = \beta\), so that \(\Phi_S = \Phi(\beta)\) is the local minimum of the potential and for all \(u \in (\ShockMin,\ShockMax)\), \(\Phi(u) > \Phi_S\), then the left hand side of \Cref{equation:equal_area_in_phi_shock_condition} is positive. By the intermediate value theorem, there must be a \(\Phi_S\) with \(\Phi(\beta)<\Phi_S<\Phi(\alpha)\), such that condition \Cref{equation:equal_area_in_phi_shock_condition} is satisfied. Consequently there is always a shock position that satisfies the equal area in \(\Phi(u)\) rule. See \Cref{fig:lemma1}.
            \end{proof}
   
            \begin{proof}[Proof of \Cref{lemma:existence_of_shock_position_equal_diffusion}: existence of continuous diffusivity shock]                
                Consider the graphs of \(D(\ShockMin)\) and \(D(\ShockMax)\) as functions of \(\Phi_S\), where  \(\Phi(\beta)<\Phi_S<\Phi(\alpha)\), \Cref{figure:continous_diffusivity_existence}. The increasing (red) curve shows the diffusivity at the upper endpoint of the shock \(\ShockMax\), while the decreasing (blue) curve shows the diffusivity at the lower endpoint of the shock \(\ShockMin\). \Cref{figure:continous_diffusivity_existence} shows the graph for the cubic diffusivity~\cref{equation:quadratic_cubic_diffusion} described in~\Cref{section:nonsymmetric_diffusion}, but is representative of any diffusivity that follows a decreasing-increasing pattern.

                \begin{figure}[htbp!]
                    \centering
                     \includegraphics[width=0.65\linewidth]{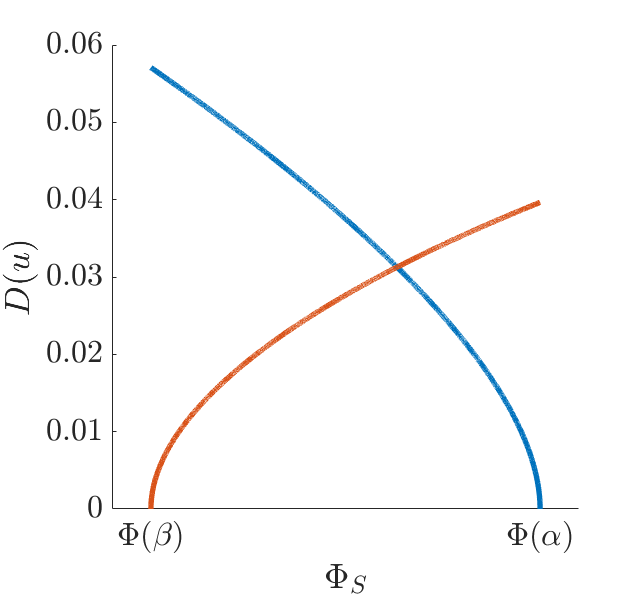}
                    \caption{Value of the cubic diffusivity~\cref{equation:quadratic_cubic_diffusion} at the shock end-points, with \(a = 0.2,\) \(b = 0.4,\) and \(\delta = 0.5\) as functions of \(\Phi_S\) for all possible shock positions. Red curve: Upper/right endpoint, \(D(\ShockMax)\). Blue curve: lower/left endpoint, \(D(\ShockMin)\).} 
                    \label{figure:continous_diffusivity_existence}
                \end{figure} 

                The smallest possible value of \(\Phi_S\) is \(\Phi(\beta)\) (see also the lowest shock in \Cref{figure:phi_shocks}). For \(\Phi_S=\Phi(\beta)\), \(\ShockMin<\alpha\) and \(\ShockMax=\beta\) so that \(D(\ShockMin)>0\) and \(D(\ShockMax)=0\). Consequently, for the lowest possible value of \(\Phi_S\), \(D(\ShockMin)>D(\ShockMax)\). The largest possible value of \(\Phi_S\) is \(\Phi(\alpha)\) (see also the uppermost shock in \Cref{figure:phi_shocks}). For \(\Phi_S=\Phi(\alpha)\), \(\ShockMin=\alpha\) and \(\ShockMax>\beta\) so that \(D(\ShockMin)=0\) and \(D(\ShockMax)>0\). Consequently, for the largest possible value of \(\Phi_S\), \(D(\ShockMin)<D(\ShockMax)\).

                Now consider the function \(D(\ShockMax) - D(\ShockMin)\) as a function of \(\Phi_S\). At \(\Phi_S = \Phi(\beta)\), \(D(\ShockMax) - D(\ShockMin)<0\), while at \(\ShockMin = \alpha\), \(D(\ShockMax) - D(\ShockMin)>0\). By the intermediate value theorem, there must be a \(\Phi_S\) where \(D(\ShockMax) - D(\ShockMin) = 0\). As a consequence, there must be at least one shock location that maintains continuity of the diffusivity~\Cref{equation:continous_diffusion_shock_condition}. If the diffusivity is monotonically decreasing for \(u<\alpha\), then \(D(\ShockMin)\) is monotonically decreasing with increasing \(\Phi_S\). Similarly, if the diffusivity is monotonically increasing for \(u>\beta\), \(D(\ShockMax)\) is monotonically increasing with increasing \(\Phi_S\). Consequently, \(D(\ShockMax) - D(\ShockMin)\) is monotonically increasing with increasing \(\Phi_S\) and so \(D(\ShockMax) - D(\ShockMin) = 0\) at only one value of \(\Phi_S\), so that the shock that maintains continuity of the diffusivity is unique.    
            \end{proof}

        \subsection{Continuous diffusivity maximises shock length} \label{proof:longest_shock}

            \Cref{lemma:longest_shock} states that the shock that ensures continuous diffusivity also maximises the shock length. .
    
            \begin{proof}[Proof of \cref{lemma:longest_shock}]
                The derivative of the shock length \(S_L\) see~\cref{equation:shock_length} with respect to \(\Phi\) is
                \begin{equation}
                    \label{equation:shock_length_derivative}%
                    S_L'(\Phi) = \ShockMax'(\Phi) - \ShockMin'(\Phi) = \frac{1}{\Phi'(\ShockMax)} - \frac{1}{\Phi'(\ShockMin)} = \frac{1}{D(\ShockMax)} - \frac{1}{D(\ShockMin)}.
                \end{equation}
                Critical points of the shock length require that \(D(\ShockMax) = D(\ShockMin) \neq 0\). To verify that continuous diffusivity maximises shock length in the decreasing-increasing case, we consider the second derivative of the shock length,
                \begin{equation}
                    \label{equation:shock_length_second_derivative}%
                    S_L''(\Phi) = \frac{D'(\ShockMin)}{D(\ShockMin)^3} - \frac{D'(\ShockMax)}{D(\ShockMax)^3}<0, \quad{\rm because}\quad D'(\ShockMin)<0\quad{\rm and}\quad D'(\ShockMax)>0.
                \end{equation}
                Since the diffusivity is continuous everywhere, the critical point is a local maximum. In the decreasing-increasing diffusivity case, there exists only one shock location with continuous diffusivity. Consequently the critical point is unique and is the global maximum.
            \end{proof}

\begin{figure}[htbp!]
                \centering
                \subcaptionbox{\label{figure:shock_positions_delta_0.5}}{\includegraphics[width=0.49\linewidth]{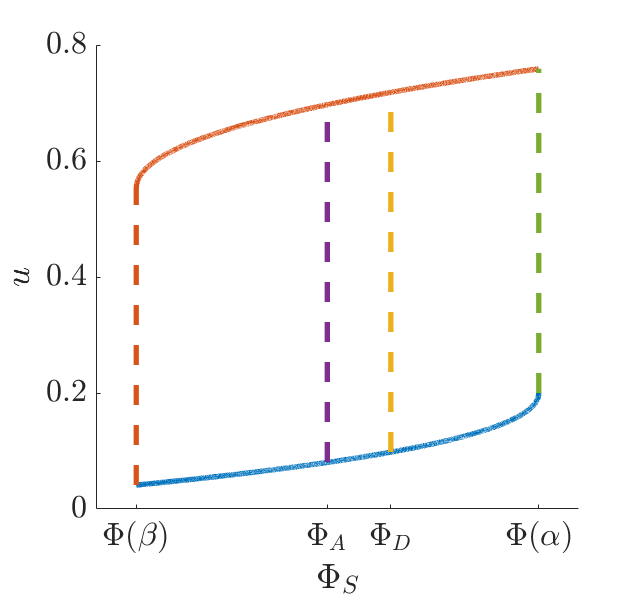}}
                \subcaptionbox{\label{figure:shock_lengths_delta_0.5}}{\includegraphics[width=0.49\linewidth]{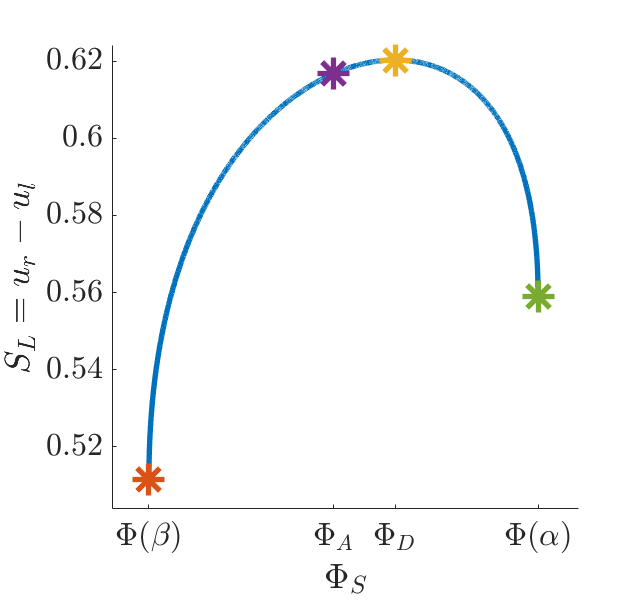}}
                \caption{Shock positions and length for the cubic diffusivity~\cref{equation:quadratic_cubic_diffusion}, with \(a = 0.2, b = 0.4\) and \(\delta = 0.5\) (a) Possible shock positions. The red curve shows possible \(\ShockMax\) values and the blue curve shows possible \(\ShockMin\) values. Possible shocks are then connections between \(\ShockMax\) and \(\ShockMin\) with constant \(\Phi(u) = \Phi_S\). (b) Shock lengths as a function of \(\Phi_S\). Red dotted line (a) and red star (b): shock where \(\ShockMax = \beta\). Purple dotted line/purple star: equal area in \(\Phi(u)\) shock. Yellow dotted line/yellow star: continuous diffusivity shock. Green dotted line/green star: shock where \(\ShockMin = \alpha\).}
            \end{figure}

\begin{figure}[htbp!]
                \centering
                \subcaptionbox{\label{figure:crazy_diffusion}}{\includegraphics[width=0.49\linewidth]{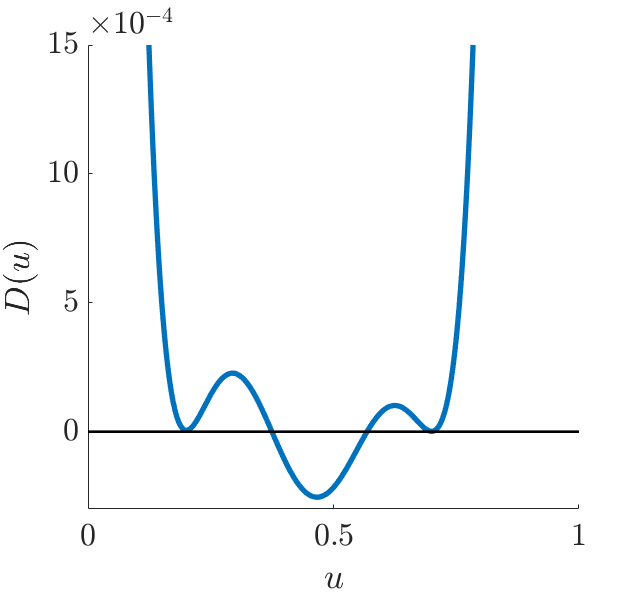}}
                \subcaptionbox{\label{figure:crazy_shock_lengths}}{\includegraphics[width=0.49\linewidth]{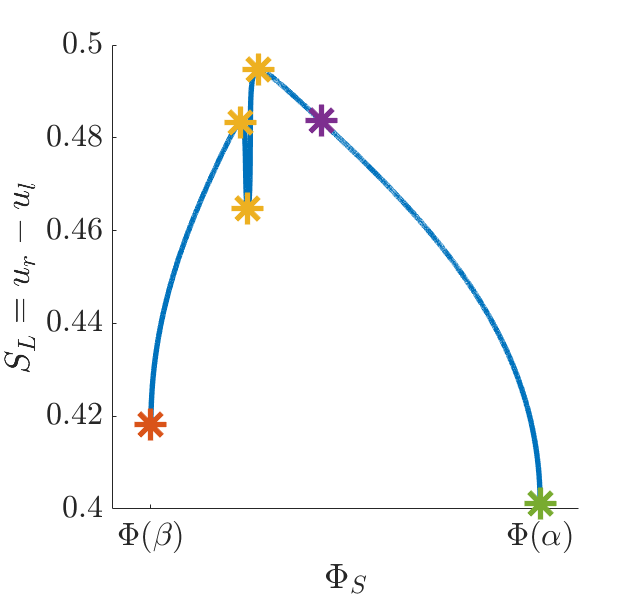}}
                \caption{Diffusivity and shocks for \(D(u)\) that is negative for \(\alpha < u < \beta,\) but does not follow the decreasing-increasing pattern. (a) Diffusivity, \(D(u).\) (b) Shock lengths as a function of \(\Phi_S\) corresponding to the diffusivity in~\Cref{figure:crazy_diffusion}. 
                Red star: shock where \(\ShockMax = \beta\). Purple star: equal are in \(\Phi(u)\) shock. Yellow stars: continuous diffusivity shock. Green star: shock where \(\ShockMin = \alpha\).}
                \label{figure:crazy}
            \end{figure}

            We illustrate~\Cref{lemma:longest_shock} using the cubic diffusivity~\Cref{equation:quadratic_cubic_diffusion} with \(\delta = 0.5\). \Cref{figure:shock_positions_delta_0.5} shows jumps in \(u\) at different values of the shock potential, \(\Phi_S\), and~\Cref{figure:shock_lengths_delta_0.5} shows shock length as a function of \(\Phi_S.\) The continuous diffusivity shock gives the longest shock. When not restricted to the decreasing-increasing profile for the diffusivity, it is possible to construct examples with multiple shock positions that satisfy a continuous diffusion rule that lie at local extrema in the shock length. \Cref{figure:crazy_diffusion} illustrates one such diffusivity that is no longer decreasing-increasing. As~\Cref{figure:crazy_shock_lengths} shows, multiple shocks satisfying continuous diffusivity are possible, and coincide with extrema of \(S_L(\Phi)\). One of these continuous diffusivity shocks is the longest possible shock.
    
        \subsection{Nonlinear regularisation selects the shock position} \label{section:Nonlinear_Regularisation_Yields_Continuous_Diffusivity_Shocks}
            \textcite{bradshaw2023geometric} used a composite regularisation to obtain shock solutions where the shock position could be moved. In this section we show that we can move shock positions using a single nonlinear regularisation term, \(-\varepsilon^2 \left(f(u)u_{xx}\right)_{xx}\). Using a travelling wave formulation, we prove \Cref{lemma:perturbation} by expanding about the shock location and \Cref{theorem:GSPT} using geometric singular perturbation theory (GSPT).

            \subsubsection{Nonlinear regularisation corresponds to a modified equal area rule}\label{section:perturbation}
                We use a perturbation argument similar to~\textcite{witelski1995shocks,witelski1996structure} to show that the nonlinear regularisation~\Cref{equation:nonlinear_fourth_order_regularisation_2_2} leads to modified equal area rule~\Cref{equation:modified_equal_area_rule} for any reaction term, \(R(u)\).

                \begin{proof}                            
                    We start by writing~\Cref{equation:nonlinear_fourth_order_regularisation_2_2} in terms of the travelling wave coordinate \(z = x - ct\),
                    \begin{equation}
                        \label{equation:nonlinear_fourth_order_regularisation_2_2_TW}%
                        -c \fd{u}{z} - \fdn{}{z}{2}\Phi(u) + \varepsilon^2 \fdn{}{z}{2}\left(f(u) \fdn{u}{z}{2}\right) - R(u) = 0.
                    \end{equation}
                    Next, to examine dynamics close to the shock we introduce a new stretched independent variable \(\xi = (z - z_s)/\varepsilon,\) where \(z_s\) is the location of the shock in \(z\). This yields
                    \begin{equation}
                        \label{equation:nonlinear_fourth_order_regularisation_2_2_Stretched_Variable_TW}%
                        -\frac{c}{\varepsilon}\fd{u}{\xi} - \frac{1}{\varepsilon^2}\fdn{}{\xi}{2}\Phi(u) + \fdn{}{\xi}{2}\left(\frac{f(u)}{\varepsilon^2}\fdn{u}{\xi}{2}\right) - R(u) = 0,    
                    \end{equation}
                    for which we consider the leading order problem as \(\varepsilon \rightarrow 0\),
                    \begin{equation}
                        \label{equation:Perturbation_leading_order_start}%
                        \fdn{}{\xi}{2}\Phi(u) - \fdn{}{\xi}{2}\left(f(u)\fdn{u}{\xi}{2}\right) = 0.        
                    \end{equation}
                    Integrating twice with respect to \(\xi,\)
                    \begin{equation}
                        \label{equation:Perturbation_leading_order_step_one}%
                        \Phi(u) - f(u)\fdn{u}{\xi}{2} = c_1 \xi + c_2.
                    \end{equation}
                    To specify the constants of integration \(c_1\) and \(c_2,\) we use the boundary conditions for \(u(\xi).\) Since \(\xi\) represents an expansion about the shock, the appropriate conditions are \(u \to u_l\) as \(\xi \to -\infty\) and \(u \to u_r\) as \(\xi \to \infty,\) with higher derivatives of \(u\) with respect to \(\xi\) vanishing as \(\xi \to \pm\infty.\) Hence, \(c_1 = 0\) and \(c_2 = \Phi_S = \Phi(\ShockMax) = \Phi(\ShockMin)\). Then,
                    \begin{equation}
                        \label{equation:Perturbation_leading_order_step_two}%
                        \fdn{u}{\xi}{2} = \frac{\Phi(u) - \Phi_S}{f(u)},
                    \end{equation}
                    with \(\Phi_S\) still to be determined. After multiplying both sides by $\d u/\d\xi$ and integrating with respect to $\xi$, we can rearrange and integrate again to obtain
                    \begin{equation}
                        \label{equation:Perturbation_leading_order_step_three}%
                        \pm \xi + c_4 = \int \frac{\df{u}}{\sqrt{2\int_{u^*}^{u} \left[\Phi(\Bar{u}) - \Phi_S\right]/f(\Bar{u}) \df{\bar{u}} + c_3}},
                    \end{equation}
                    This is similar to the condition that~\textcite{pego1989front} showed yields an equal area rule in \(\Phi(u)\). For \(u = \ShockMax\) and \(u = \ShockMin\), we want \(\xi \rightarrow \pm \infty\) so we require the outer integral to diverge for \(u = \ShockMax\) and \(u = \ShockMin\). Thus, the square root term on the denominator should go to zero at \(u = \ShockMax\) and \(u = \ShockMin\), so
                    \begin{equation}
                        2\int_{u^*}^{\ShockMin} \left(\frac{\Phi(\Bar{u}) - \Phi_S}{f(\Bar{u})}\right)\df{\Bar{u}} + c_3 = 0 \qquad \text{and} \qquad 2\int_{u^*}^{\ShockMax} \left(\frac{\Phi(\Bar{u}) - \Phi_S}{f(\Bar{u})}\right)\df{\Bar{u}} + c_3 = 0,
                    \end{equation}
                    which implies
                    \begin{equation}                    
                        \int_{\ShockMin}^{\ShockMax} \left(\frac{\Phi(\Bar{u}) - \Phi_S}{f(\Bar{u})}\right)\df{\Bar{u}} = 0.
                    \end{equation}
                    This equation~\cref{equation:modified_equal_area_rule} is the modified equal area rule.
                \end{proof}
        
            \subsubsection{Existence and uniqueness of travelling wave solutions with continuous diffusivity shocks for a cubic reaction term}\label{section:GSPT}
        
                \begin{proof}[Proof of \Cref{theorem:GSPT} using geometric singular perturbation theory]
                
                    Following a similar argument to~\textcite{li2021shock}, we use geometric singular perturbation theory (GSPT)~\parencite{fenichel1979geometric,jones1995geometric,kuehn2015multiple,bradshaw2023geometric,Harley2014} to prove the existence of travelling wave solutions satisfying the modified equal area rule~\Cref{equation:modified_equal_area_rule} in the \(\varepsilon \to 0\) limit of the reaction--diffusion equation with nonlinear regularisation~\Cref{equation:nonlinear_fourth_order_regularisation_2_2}. Writing~\Cref{equation:nonlinear_fourth_order_regularisation_2_2} in the travelling coordinate, \(z=x-ct\), we obtain
                    \begin{equation}
                        \label{eq:TW_GSPT}%
                        \fd{}{z}\left[-c u - \fd{}{z}\Phi(u) + \varepsilon^2 \fd{}{z}\left(f(u) \fdn{u}{z}{2}\right)\right] - R(u) = 0.
                    \end{equation}
                    To write the fourth-order equation~\cref{eq:TW_GSPT} as a system of four first-order ordinary differential equations, we define the new variables
                    \begin{subequations}
                        \begin{align}
                            p &:= -c u - \fd{}{z}\Phi(u) + \varepsilon^2 \fd{}{z}\left[f(u) \fdn{u}{z}{2}\right], \label{equation:p} \\
                            v &:= - \Phi(u) + \varepsilon^2 f(u) \fdn{u}{z}{2}, \label{equation:v} \\
                            w &:= \varepsilon \fd{u}{z}. \label{equation:w}
                        \end{align}
                    \end{subequations}
                    We then obtain the four-dimensional dynamical system with two slow variables, \(p\) and \(v\), and two fast variables, \(u\) and \(w\),
                    \begin{subequations}
                        \label{equations:slow_system}%
                        \begin{align}
                            \varepsilon\fd{u}{z} &= w, \label{equations:slow_system_u} \\
                            \varepsilon \fd{w}{z} &= \frac{v + \Phi(u)}{f(u)}, \label{equations:slow_system_w}\\
                            \fd{p}{z} &= R(u), \label{equations:slow_system_p} \\
                            \fd{v}{z} &= p + cu. \label{equations:slow_system_v}
                        \end{align}
                    \end{subequations}
                    When \(f(u) = 1\), this is identically the system for the regularised reaction--diffusion equation~\Cref{equation:standard_fourth_order_regularisation} studied by~\textcite{li2021shock}. Since the small parameter \(\varepsilon\) multiplies the derivative in~\cref{equations:slow_system_u,equations:slow_system_w}, the system~\cref{equations:slow_system} is singularly perturbed. We refer to~\cref{equations:slow_system} as the slow system. Like in~\Cref{section:perturbation}, we can obtain a topologically-equivalent system to~\cref{equations:slow_system} by introducing the fast variable \(\xi = z/\varepsilon\), where without loss of generality we set the shock position \(z_s = 0\). Then, we obtain the fast system,
                    \begin{subequations}
                        \label{equations:fast_system}%
                        \begin{align}
                            \fd{u}{\xi} &= w, \\
                            \fd{w}{\xi} &= \frac{v + \Phi(u)}{f(u)}, \\
                            \fd{p}{\xi} &= \varepsilon R(u), \\
                            \fd{v}{\xi} &= \varepsilon(p + cu).
                        \end{align}
                    \end{subequations}
                    For a cubic reaction term,
                    \begin{equation}
                        \label{equation:cubic_reaction}%
                        R(u) = u(1 - u)\left(u - \gamma\right),
                    \end{equation}
                    the three critical points for the equivalent systems~\cref{equations:fast_system,equations:slow_system} in four dimensional \((u,w,p,v)\) space are
                    \begin{subequations}
                        \begin{align}
                            u_0 &= (0,0,0,-\Phi(0)), \label{equation:critical_point_0} \\
                            u_\gamma &= (\gamma,0,-c\gamma,-\Phi(\gamma)), \label{equation:critical_point_gamma}\\
                            u_1 &= (1,0,-c,-\Phi(1)).\label{equation:critical_point_1}
                        \end{align}
                    \end{subequations}
                    With \(f(0) > 0\) and \(f(1) > 0\) the critical points \(u_1\) and \(u_0\) are both saddles, since in both cases the Jacobian has two positive eigenvalues and two negative eigenvalues. A travelling wave solution to the reaction--diffusion equation corresponds to a heteroclinic connection between the equilibria \(u_1\) and \(u_0\), which represent end-states of the wave as \(z \to -\infty\) and \(z \to \infty\) respectively. We use GSPT to establish existence of these travelling wave solutions.
    
                    To apply GSPT, we first take the singular \(\varepsilon = 0\) limit in the slow~\cref{equations:slow_system} and fast~\cref{equations:fast_system} systems. As \(\varepsilon \to 0\) the four-dimensional slow system~\cref{equations:slow_system} yields the two-dimensional reduced problem with two algebraic constraints,
                    \begin{subequations}
                        \label{equations:reduced_problem_full}%
                        \begin{align}
                            0 &= w, \label{equations:reduced_problem_full_u}\\
                            0 &= \frac{v + \Phi(u)}{f(u)}, \label{equations:reduced_problem_full_w}\\
                            \fd{p}{z} &= R(u), \label{equations:reduced_problem_full_p} \\
                            \fd{v}{z} &= p + cu. \label{equations:reduced_problem_full_v}
                        \end{align}
                    \end{subequations}
                    In the same \(\varepsilon \to 0\) limit, the fast subsystem becomes a two-dimensional layer problem with two parameters,
                    \begin{subequations}
                        \label{equations:layer_problem_full}%
                        \begin{align}
                            \fd{u}{\xi} &= w, \label{equations:layer_problem_full_u}\\
                            \fd{w}{\xi} &= \frac{v + \Phi(u)}{f(u)}, \label{equations:layer_problem_full_w}\\
                            \fd{p}{\xi} &= 0, \label{equations:layer_problem_full_p} \\
                            \fd{v}{\xi} &= 0. \label{equations:layer_problem_full_v}
                        \end{align}
                    \end{subequations}
                    We then analyse the layer and reduced problems independently to obtain singular limit solutions. If these exist, since in the regularisation \(0 < \varepsilon \ll 1\), Fenichel theory guarantees that these singular limit solutions persist as solutions to the regularised reaction--diffusion equation~\cref{equation:nonlinear_fourth_order_regularisation_2_2}.
    
                    \paragraph{The layer problem}
                        We begin by considering the layer problem~\cref{equations:layer_problem_full}. Since the fast variable \(\xi\) stretches the solution close to the shock, the layer problem describes the dynamics inside the shock. Flow in the layer problem is independent of the slow variables \(p\) and \(v\), which equations~\cref{equations:layer_problem_full_p,equations:layer_problem_full_v} show are constants of integration. The dynamics are then governed by the two-dimensional system for the fast variables,
                        \begin{subequations}
                            \label{equations:layer_problem}%
                            \begin{align}
                                \fd{u}{\xi} &= w, \\
                                \fd{w}{\xi} &= \frac{v + \Phi(u)}{f(u)}. 
                            \end{align}
                        \end{subequations}
                        Throughout the analysis based on GSPT, an important entity is the critical manifold corresponding to the equilirbia of the layer problem,
                        \begin{equation}
                            \label{equation:critical_manifold}%
                            \mathcal{S} = \{(u,w,p,v) \in \mathbb{R}^4 : w = 0, v = - \Phi(u)\}.
                        \end{equation}
                        In 4D, layer flow along lines of constant \(p\) and \(v\) (along fast fibres) terminates on \(\mathcal{S}.\) To explore the flow along fast fibres near the critical manifold, we consider the Jacobian of the layer problem~\cref{equations:layer_problem},
                        \begin{equation}
                            \label{equation:layer_problem_jacobian}%
                            J = \begin{bmatrix} 0 & 1 \\ \left(D(u)f(u) - (v + \Phi(u))f'(u)\right)/f(u)^2 & 0 \end{bmatrix}. 
                        \end{equation}
                        Evaluating \(J\) at equilibria of the layer problem on the critical manifold (where \(v = - \Phi(u)\)), we obtain the eigenvalues
                        \begin{equation}
                            \label{equation:layer_problem_eigenvalue}%
                            \lambda_{\pm} = \pm \sqrt{\frac{D(u)}{f(u)}}. 
                        \end{equation}
                        If \(f(u) > 0\) for \(u \in [0,1]\), the eigenvalues~\cref{equation:layer_problem_eigenvalue} remain real when \(D(u) \geq 0\) and become purely imaginary when \(D(u) < 0\). Consequently, the critical manifold is normally hyperbolic in two regions where diffusivity is positive. We define these regions to be \(\mathcal{S}_0\) and \(\mathcal{S}_1\), such that \(\mathcal{S}_0\) refers to regions of \(\mathcal{S}\) where \(0 < u < \alpha,\) and \(\mathcal{S}_1\) refers to regions of \(\mathcal{S}\) where \(\beta < u < 1.\) In contrast, the critical manifold is not normally hyperbolic in the region \(\mathcal{S}_{neg}\) where \(D(u) < 0,\) that is for \(u \in (\alpha,\beta).\) Since the shock jumps over the region of negative diffusivity, equilibria of the layer problem lie on normally hyperbolic regions of the critical manifold. A solution to the layer problem~\cref{equations:layer_problem} is then a heteroclinic connection between a point on \(\mathcal{S}_0\) and a point on \(\mathcal{S}_1,\) which describes the shock that jumps over \(\mathcal{S}_{neg}.\) 
                    
                        We obtain the heteroclinic connection that solves the layer problem by exploiting the Hamiltonian structure of the layer problem~\cref{equations:layer_problem}. The Hamiltonian is
                        \begin{equation}
                            \label{equation:layer_problem_hamiltonian}%
                            H(u,w) = -\frac{1}{2}w^2 v F(u) + G(u),
                        \end{equation}
                        where \(F(u)\) and \(G(u)\) are the antiderivatives of \(1/f(u)\) and \(\Phi(u)/f(u)\) respectively. 
                        A solution to the layer problem is confined to a level set of the Hamiltonian and hence requires
                        \begin{equation}
                            -\frac{1}{2}w^2 + v F(u_l) + G(u_l) = -\frac{1}{2}w^2 + v F(u_r) + G(u_r).
                        \end{equation}
                        Given the definition of $F(u)$ and $G(u)$, this is equivalent to 
                        \begin{equation}
                            \int_{u_l}^{u_r} \frac{\Phi(u) + v}{f(u)} \df{u} = 0.
                        \end{equation}
                        Since \(v = -\Phi(\ShockMin) = -\Phi(\ShockMax) = -\Phi_S\) this is the modified equal area rule, equation \eqref{equation:modified_equal_area_rule}.
                                 
                        With \(f(u) > 0\) for \(u \in [0,1]\) the modified equal area rule cannot give shocks that jump to or from the zeros of the diffusivity, \(\alpha,\) \(\beta\). For these shocks, \(\Phi(u) - \Phi_S\) is either always positive or always negative so that \([\Phi(u) - \Phi_S]/f(u)\) remains either always positive or always negative. Consequently, the integral can never be zero and these shock positions cannot satisfy an equal area rule. This guarantees that \(\ShockMin < \alpha\) and \(\ShockMax > \beta\) will hold for all shocks from the modified equal area rule. Away from the shock, the solution will remain on \(\mathcal{S}_0\) and \(\mathcal{S}_1,\) the normally hyperbolic parts of the critical manifold. These dynamics away from the shock are governed by the reduced problem, which we now consider.
          
                    \paragraph{The reduced problem}
                        To understand dynamics outside the shock, we consider the two-dimensional form of the reduced problem,
                        \begin{subequations}
                            \label{equations:reduced_problem}%
                            \begin{align}            
                                \fd{p}{z} &= R(u), \\
                                \fd{v}{z} &= p + cu,
                            \end{align}
                        \end{subequations}
                        which according to the full reduced problem~\cref{equations:reduced_problem_full} must evolve confined to the critical manifold \(\mathcal{S}.\) To analyse the reduced problem, we use the approach of~\textcite{li2020travelling} and rewrite the reduced problem~\cref{equations:reduced_problem} as
                            \begin{subequations}
                                \label{equations:reduced_problem_two}%
                                \begin{align}
                                    \fd{p}{z} &= R(u), \\
                                    -D(u) \fd{u}{z} &= p + cu.
                                \end{align}
                            \end{subequations}    
                        Since this form introduces singularities at \(u = \alpha\) and \(u = \beta\) where \(D(u) = 0,\) we desingularise by introducing a new variable \(\psi\) such that \(\d\psi = \d z/D(u)\). This yields
                        \begin{subequations}
                            \label{equations:phase_plane}%
                            \begin{align}
                                \fd{p}{\psi} &= R(u)D(u), \\
                                \fd{u}{\psi} &= -p - cu,
                            \end{align}
                        \end{subequations}
                        from which we investigate solutions in the \((p,u)\) phase plane.
    
                        A solution to the reduced problem connects the unstable manifold of the equilibrium \(u_1\) to the stable manifold of \(u_0\). To apply Fenichel theory, flow in the reduced problem must remain on \(\mathcal{S}_0\) or \(\mathcal{S}_1\), the normally hyperbolic regions of \(\mathcal{S}\). In the reduced problem flow that originates at a point on either \(\mathcal{S}_0\) or \(\mathcal{S}_1\) must remain on that part of the critical manifold. Thus, in general a single trajectory emanating from \(u_1\) is unable to connect directly to \(u_0\). Therefore, to connect flow on \(\mathcal{S}_0\) to flow on \(\mathcal{S}_1\), we require a shock that jumps across \(\mathcal{S}_{neg},\) the region where \(\mathcal{S}\) is not normally hyperbolic. This is provided by the layer problem, which connects regions of \(\mathcal{S}_0\) or \(\mathcal{S}_1\) along fast fibres of constant \(p\) and \(v\). As shown when analysing the layer problem, these shocks satisfy the modified equal area rule~\cref{equation:modified_equal_area_rule}.
    
                        The phase planes in~\Cref{figures:phase_planes} illustrate the concept of adding a shock to trajectories of the reduced problem. In each phase plane of~\Cref{figures:phase_planes}, the red vertical lines indicate the zeros of the diffusivity, and the continuous diffusivity shock must occur at the black dashed vertical lines. In all scenarios, trajectories along the unstable manifold of \(u_1 = (1,-c)\) do not connect with trajectories entering the stable manifold of \(u_0 = (0,0).\) The continuous diffusivity shock is possible when the two trajectories are such that the fast variable \(p\) takes the same value when  \(u = \ShockMin\) and \(u = \ShockMax,\) where \(\ShockMin\) and \(\ShockMax\) take the fixed values corresponding to the continuous diffusivity shock. The shock is then the horizontal line with constant \(p\) that connects the unstable manifold of \(u_1\) at \(u = \ShockMax\) to the stable manifold of \(u_0\) at \(u = \ShockMin.\) This shock is shown in~\Cref{figure:phase_plane_equal_diffusion_right_speed}. Furthermore, there is a unique travelling wave speed, \(c,\) that is necessary for the continuous diffusivity shock to be possible.~\Cref{figure:phase_plane_equal_diffusion_too_slow,figure:phase_plane_equal_diffusion_too_fast} illustrate two phase planes with the incorrect speed \(c,\) such that the line connecting the two trajectories at \(u = \ShockMax\) and \(u = \ShockMin\) is not horizontal, and \(p\) is not conserved across the connection. Since non-constant \(p\) contradicts the layer flow, the connection is not a permissible shock. This illustrates that the nonlinear regularisation and continuous diffusivity shock rule not only selects the shock position, but also the speed of advance for the travelling wave solution.                    
    
                        \begin{figure}[htbp!]        
                            \centering            
                            \subcaptionbox{\label{figure:phase_plane_equal_diffusion_too_slow}}{\includegraphics[width=0.49\linewidth]{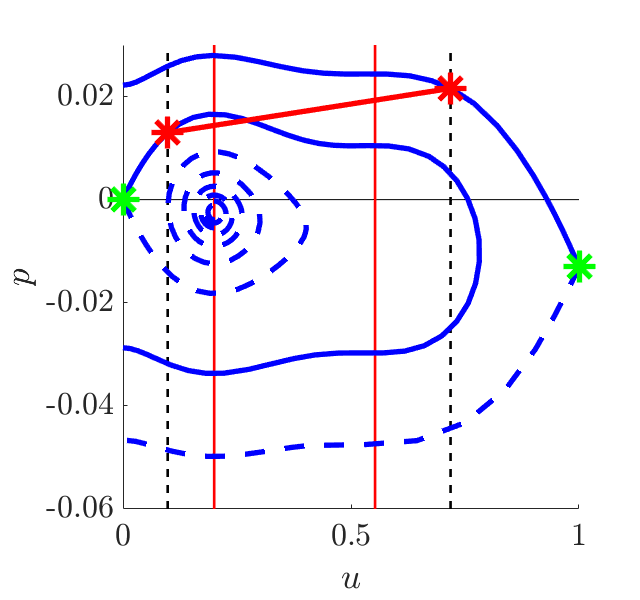}}
                            \subcaptionbox{\label{figure:phase_plane_equal_diffusion_too_fast}}{\includegraphics[width=0.49\linewidth]{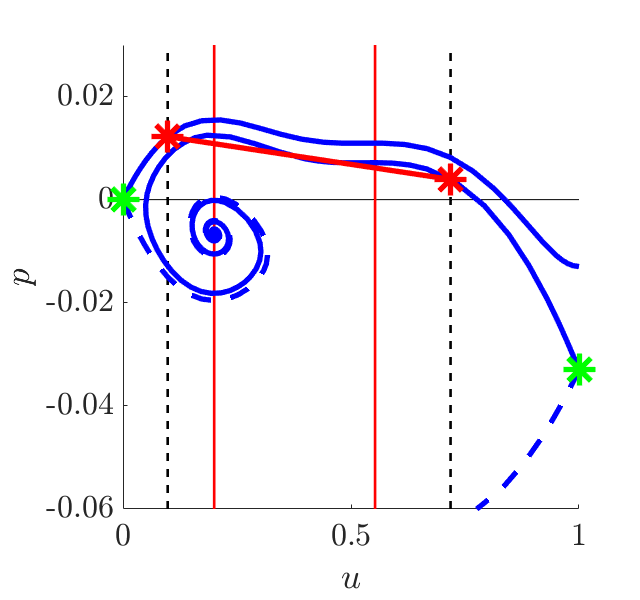}}
                            \subcaptionbox{\label{figure:phase_plane_equal_diffusion_right_speed}}{\includegraphics[width=0.49\linewidth]{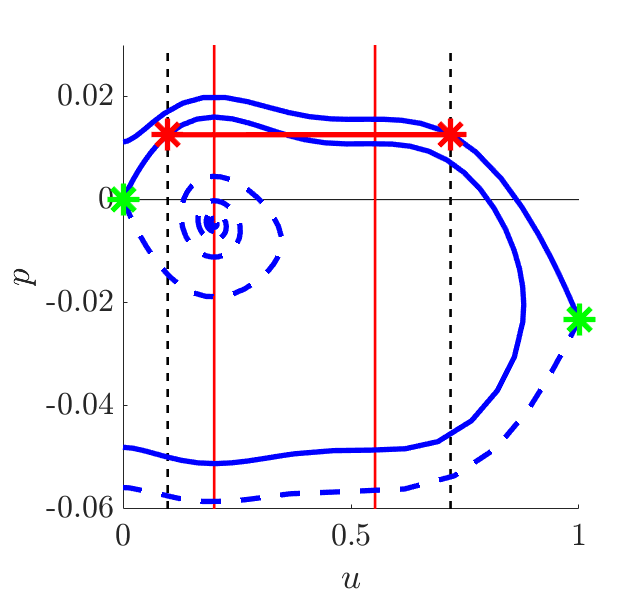}}
                            \subcaptionbox{\label{figure:phase_plane_speed}}{\includegraphics[width=0.49\linewidth]{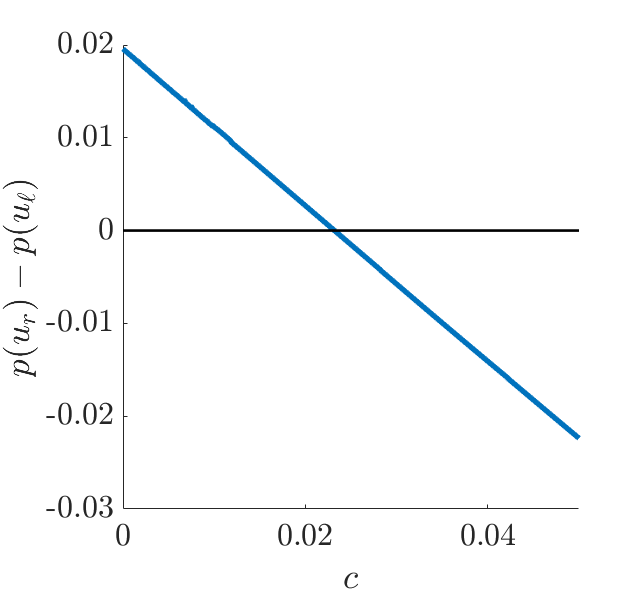}}
                            \caption{Phase planes of the reduced problem~\Cref{equations:phase_plane}. The green stars are at \(u_1 = (1,-c)\) and \(u_0 = (0,0)\). Solid blue curves are the unstable manifold of \(u_1\) and stable manifold of \(u_0,\) trajectories corresponding to the travelling wave solution. Dashed trajectories are the stable manifold of \(u_1\) and unstable manifold of \(u_0,\) which the travelling wave does not access. The vertical red lines are the zeros of the diffusivity, and the black dashed vertical lines show the location of \(\ShockMax\) and \(\ShockMin\) for the continuous diffusivity shock. The red stars show where the unstable manifold of \(u_1\) crosses the \(\ShockMax\) line and where the stable manifold of \(u_0\) crosses the \(\ShockMin\) line. The red line connecting these end points is a jump corresponding to a potential shock, but only satisfies the layer dynamics if \(p\) is conserved across the shock. (a) \(c = 0.013\), for which we do not obtain a valid continuous diffusivity shock. (b) \(c = 0.033\), for which we do not obtain a valid continuous diffusivity shock. (c) \(c = 0.0232\), for which the red line is horizontal, meaning the continuous diffusivity shock corresponds to this speed. (d) Wave speed \(c\) versus \(p(u_r) - p(u_\ell).\) The correct speed occurs when \(p(u_r) - p(u_\ell) = 0\).}
                            \label{figures:phase_planes}
                        \end{figure}
    
                        We also mention that these travelling waves are weak solutions to the original reaction--diffusion equation,~\cref{equation:reaction_diffusion_equation}. Considering weak solutions to reaction--diffusion equations following~\parencite{whitham2011linear,witelski1995shocks} we find they must conserve \(\Phi(u)\) across the shock and we also require
                        \begin{equation}
                            \label{equation:second_weak_condition}%
                            c = -\frac{\Phi(\ShockMax)_z - \Phi(\ShockMin)_z}{\ShockMax - \ShockMin} = -\frac{1}{\ShockMax - \ShockMin}\left[D(\ShockMax)\fd{u}{z}\Big|_{\ShockMax} - D(\ShockMin) \fd{u}{z}\Big|_{\ShockMin}\right],
                        \end{equation}
                        which can be considered an analogy for parabolic equations to the Rankine--Hugoniot condition for hyperbolic equations. From the reduced problem~\Cref{equations:reduced_problem_two}, \(\Phi(u)_z = D(u)u_z = - p - cu\). Since \(p\) is conserved across the shock, continuous diffusivity shock solutions satisfy the weak solution condition~\Cref{equation:second_weak_condition}. Hence, given that we have established heteroclinic connections separately in the layer and reduced problems, GSPT establishes the existence of the shock-fronted travelling wave solutions to the reaction--diffusion equation with nonlinear regularisation,~\Cref{equation:nonlinear_fourth_order_regularisation_2_2}.
                \end{proof}

            \subsubsection{Regularisation with exponential decay is an example that yields continuous diffusivity shocks}\label{section:Exponential_Regularisation}
\begin{figure}[htbp!]
                    \centering
                    \subcaptionbox{}{\includegraphics[width = 0.51\textwidth]{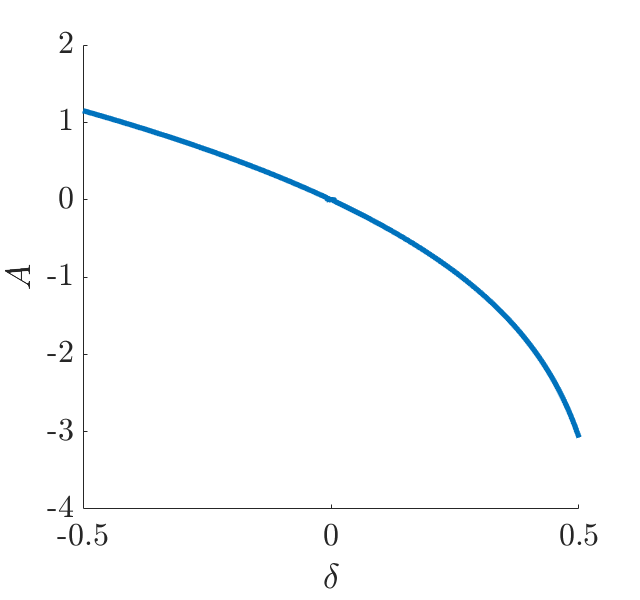}}
                    \subcaptionbox{\label{figure:exponential_f}}{\includegraphics[width=0.49\linewidth]{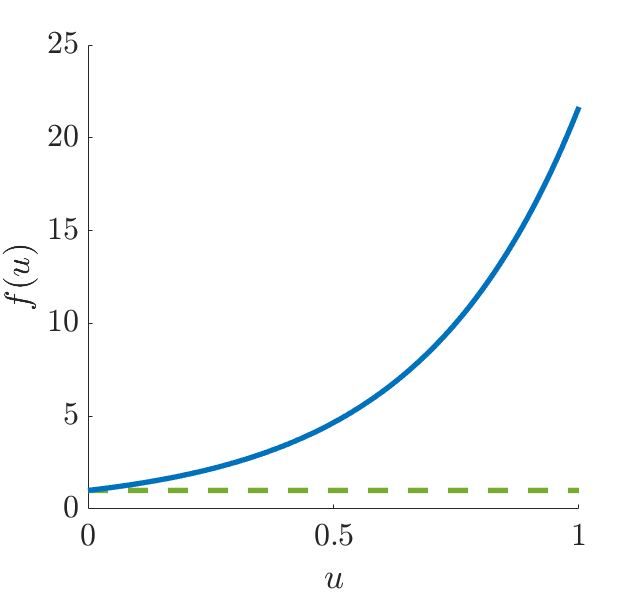}}
                    \subcaptionbox{\label{figure:modifiedEqualArea}}{\includegraphics[width=0.49\linewidth]{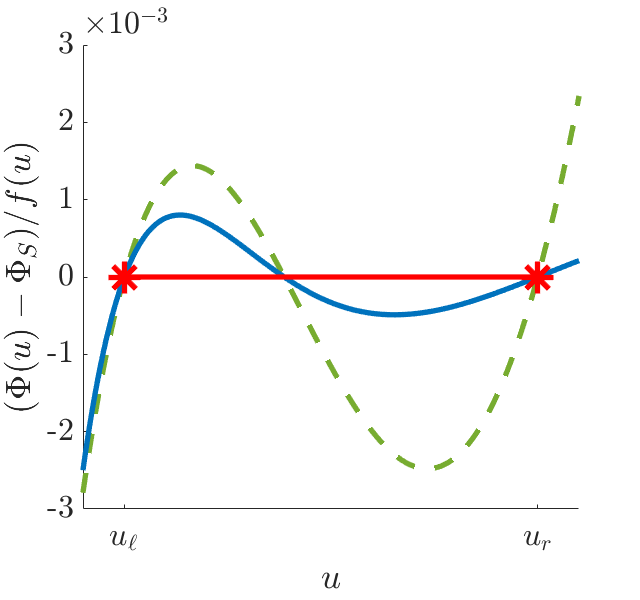}}
                    \caption{The modified equal area rule~\cref{equation:modified_equal_area_rule} for quadratic/cubic diffusivity ~\Cref{equation:quadratic_cubic_diffusion} with \(a = 0.2\), \(b = 0.4,\) and exponential regularisation \(f(u) = \e^{-Au}.\) (a) Values of \(A\) that give a continuous diffusivity shock versus \(\delta\) (b) The regularisation function \(f(u).\) Green dotted horizontal line: \(f(u) = 1\). Blue curve: \(f(u) = \e^{-Au}\) with \(A = -3.0757\). (c) The modified equal area rule. Red stars represent \(\ShockMin\) and \(\ShockMax\) for the continuous diffusivity shock. The green dotted curve shows \((\Phi(u) - \Phi_S)/f(u)\) where \(f(u) = 1\) and \(\Phi_S = \Phi(\ShockMin) = \Phi(\ShockMax)\). For \(f(u) = 1,\) the modified equal area rule~\Cref{equation:modified_equal_area_rule} reduces to the standard equal area in \(\Phi(u)\) rule \Cref{equation:equal_area_in_phi_shock_condition}, and is not satisfied by this shock as the area in the lobes of the green curves cut off by the red horizontal line are not equal. The blue curve is \((\Phi(u) - \Phi_S)/f(u)\) with \(f(u) = \e^{-Au}\) and \(A = 3.0757\). This has rescaled the areas in the lobes such that the areas under of the lobes in the blue curve are equal and so satisfies the modified equal area rule.}
                    \label{figure:A_values}
                \end{figure}

                We now demonstrate that the modified equal area rule~\cref{equation:modified_equal_area_rule} can be used to yield a travelling wave solution with a continuous diffusivity shock. As an example, we use the exponential function
                \begin{equation}
                    \label{equation:exponential_f}%
                    f(u) = \e^{-Au},        
                \end{equation}
                which satisfies the \(f(u) > 0\) positivity requirement. For a polynomial diffusivity \(D(u)\), we can also compute the integral for the modified equal area rule explicitly with the exponential function,
                \begin{equation}
                \label{equation:modified_equal_area_rule_integrated}%
                    \int_{\ShockMin}^{\ShockMax} \frac{\left(\Phi(u) - \Phi_S\right)}{\e^{-Au}}\df{u} = \frac{\e^{A\ShockMax}}{A} \sum_{i = 1}^{n + 1} \left(-1\right)^i\frac{\Phi(\ShockMax)^{(i)}}{A^i} - \frac{\e^{A\ShockMin}}{A}\sum_{i = 1}^{n + 1} \left(-1\right)^i\frac{\Phi(\ShockMin)^{(i)}}{A^i} = 0,
                \end{equation}
                where \(n\) is the degree of \(D(u)\) and \(\Phi(u)^{(i)}\) means the \(i\)-th derivative of \(\Phi(u)\) with respect to \(u\). So given the end-points \(\ShockMin\) and \(\ShockMax\) of the continuous diffusivity shock, to obtain the correct regularisation we need to find a value of \(A\) which satisfies \Cref{equation:modified_equal_area_rule_integrated}. \Cref{figure:A_values} shows values of \(A\) corresponding to a continuous diffusivity rule found numerically for the cubic diffusivity~\Cref{equation:quadratic_cubic_diffusion} where \(a = 0.2\) and \(b = 0.4\). For \(\delta = 0.5\) we find that \(A\approx 3.0757\) will give a shock position with continuous diffusivity. \Cref{figure:exponential_f} shows this exponential, and~\Cref{figure:modifiedEqualArea} visualises how the modified equal area rule is satisfied. For each \(\delta\) we can find an \(A\) that satisfies the modified equal area rule~\Cref{equation:modified_equal_area_rule}, as the following proof of~\Cref{lemma:finding_A} demonstrates.
            
                \begin{proof}
                    To prove that we can always find a value of \(A\) satisfying the modified equal area rule, we define
                    \begin{equation}
                    \label{eq:G}%
                        G(A) = \int_{\ShockMin}^{\ShockMax} \frac{\left(\Phi(u) - \Phi_S\right)}{\e^{-Au}} \df{u},
                    \end{equation}
                    and seek a value of \(A\) such that \(G(A) = 0.\) For \(A = 0,\) \(\e^{-Au} = 1\) so \(G(0) = \int_{\ShockMin}^{\ShockMax} \left(\Phi(u) - \Phi_S\right) \df{u},\) which is a constant. When \(A \neq 0,\)
                    \begin{equation}
                        \label{eq:G_evaluated}%
                        G(A) = \frac{\e^{A\ShockMax}}{A}\left(\sum_{i = 1}^{n + 1} \left(-1\right)^i\frac{\Phi(\ShockMax)^{(i)}}{A^i}\right) - \frac{\e^{A\ShockMin}}{A}\left(\sum_{i = 1}^{n + 1} \left(-1\right)^i\frac{\Phi(\ShockMin)^{(i)}}{A^i}\right).
                    \end{equation}
                    Taking the limit as \(A \to \infty\) gives
                    \begin{align*}
                        \lim_{A\rightarrow\infty} G(A) &= \lim_{A\rightarrow\infty}\left(\frac{\e^{A\ShockMax}}{A}\left(\sum_{i = 1}^{n + 1} \left(-1\right)^i\frac{\Phi(\ShockMax)^{(i)}}{A^i}\right) - \frac{\e^{A\ShockMin}}{A}\left(\sum_{i = 1}^{n + 1} \left(-1\right)^i\frac{\Phi(\ShockMin)^{(i)}}{A^i}\right)\right) \\
                        &= \lim_{A\rightarrow\infty} \e^{A\ShockMax}\left(\frac{1}{A}\left(\sum_{i = 1}^{n + 1} \left(-1\right)^i\frac{\Phi(\ShockMax)^{(i)}}{A^i}\right) - \frac{\e^{A(\ShockMin -\ShockMax)}}{A}\left(\sum_{i = 1}^{n + 1} \left(-1\right)^i\frac{\Phi(\ShockMin)^{(i)}}{A^i}\right)\right).
                    \end{align*}
                    To leading order,
                    \begin{equation}
                        \lim_{A\rightarrow\infty} G(A) = \lim_{A\rightarrow\infty} \e^{A\ShockMax}\left(-\frac{\Phi'(\ShockMax)}{A^2}\right) = \lim_{A\rightarrow\infty} \e^{A\ShockMax}\left(-\frac{D(\ShockMax)}{A^2}\right) = -\infty,
                    \end{equation}  
                    because \(D(\ShockMax) > 0,\) and so the term inside the brackets is negative. We now consider the limit as \(A \rightarrow -\infty,\) which yields
                    \begin{align*}
                        \lim_{A\rightarrow-\infty} G(A) &= \lim_{A\rightarrow-\infty}\left(\frac{\e^{A\ShockMax}}{A}\left(\sum_{i = 1}^{n + 1} \left(-1\right)^i\frac{\Phi(\ShockMax)^{(i)}}{A^i}\right) - \frac{\e^{A\ShockMin}}{A}\left(\sum_{i = 1}^{n + 1} \left(-1\right)^i\frac{\Phi(\ShockMin)^{(i)}}{A^i}\right)\right) \\
                        &= \lim_{A\rightarrow-\infty} \e^{A\ShockMin}\left(\frac{\e^{A(\ShockMax - \ShockMin)}}{A}\left(\sum_{i = 1}^{n + 1} \left(-1\right)^i\frac{\Phi(\ShockMax)^{(i)}}{A^i}\right) - \frac{1}{A}\left(\sum_{i = 1}^{n + 1} \left(-1\right)^i\frac{\Phi(\ShockMin)^{(i)}}{A^i}\right)\right).
                    \end{align*}
                    To leading order,
                    \begin{equation}
                        \lim_{A\rightarrow-\infty} G(A) = \lim_{A\rightarrow-\infty} \e^{A\ShockMin}\left(\frac{\Phi'(\ShockMin)}{A^2}\right) = \lim_{A\rightarrow-\infty} e^{A\ShockMin}\left(\frac{D(\ShockMin)}{A^2}\right) \rightarrow 0^+,
                    \end{equation}
                    because \(D(\ShockMin) > 0.\) Since \(G(A) \to -\infty\) as \(A \to -\infty\), \(G(A) \rightarrow 0^+\) as \(A \to \infty\), and \(G(A)\) is continuous, the intermediate value theorem guarantees that there exists a value of \(A\) such that \(G(A) = 0.\)
                \end{proof}

\begin{figure}[htbp!]
                    \centering
                    \subcaptionbox{\label{figure:quadratic_f}}{\includegraphics[width=0.49\linewidth]{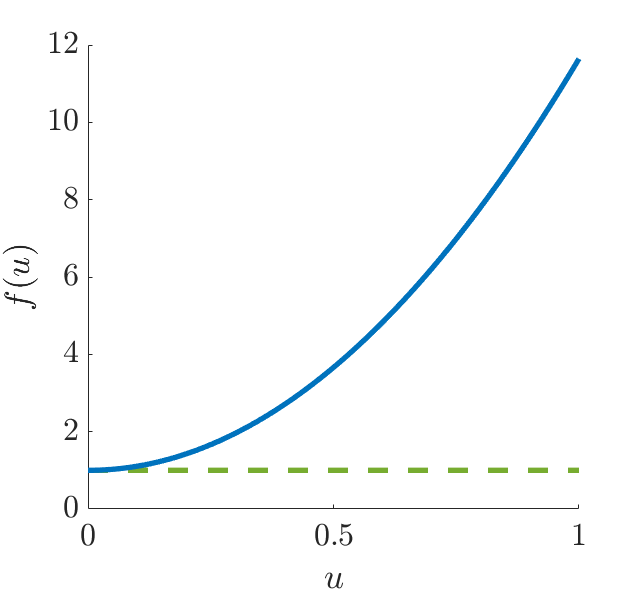}}
                    \subcaptionbox{\label{figure:quadraticModifiedEqualArea}}{\includegraphics[width=0.49\linewidth]{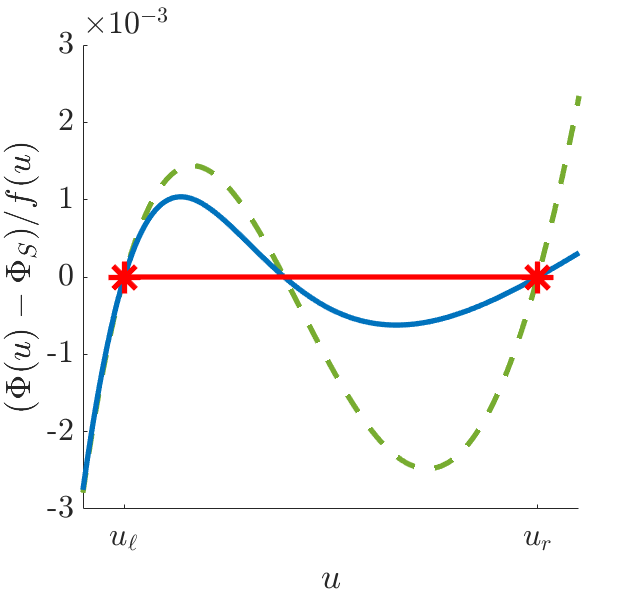}}
                    \caption{The modified equal area rule~\cref{equation:modified_equal_area_rule} for quadratic/cubic diffusivity ~\Cref{equation:quadratic_cubic_diffusion} with \(a = 0.2\), \(b = 0.4,\) and quadratic regularisation \(f(u) = 1 + Au^2\).
                    (a) The regularisation function, \(f(u).\) Green dotted horizontal line: \(f(u) = 1\). Blue curve: \(f(u) = 1 + Au^2\) with \(A = 10.6453\). (b) The modified equal area rule. Red stars represent \(\ShockMin\) and \(\ShockMax\) for the continuous diffusivity shock. The green dotted curve shows \((\Phi(u) - \Phi_S)/f(u)\) where \(f(u) = 1\) and \(\Phi_S = \Phi(\ShockMin) = \Phi(\ShockMax)\). For \(f(u) = 1,\) the standard equal area in \(\Phi(u)\) rule \Cref{equation:equal_area_in_phi_shock_condition} is not satisfied. The blue curve is \((\Phi(u) - \Phi_S)/f(u)\) with \(f(u) = 1 + Au^2\) and \(A = 10.6453\), which satisfies the modified equal area rule.}
                    \label{figures:modifiedEqualArea_all}
                \end{figure}

                Although~\Cref{lemma:finding_A} is for the exponential function, any other positive function \(f(u)\) can also be used in the modified equal area rule. For example, the quadratic function
                \begin{equation}
                    \label{equation:quadratic_f}%
                    f(u) = 1 + A u^2, \quad \textrm{ with } \quad A \approx 10.6453,
                \end{equation}
                is another possibility, as illustrated in~\Cref{figure:quadratic_f,figure:quadraticModifiedEqualArea}.
                The value of \(A\) can also be adjusted to obtain a continuum of shock positions like the composite regularisation of~\textcite{bradshaw2023geometric}, not just the continuous diffusivity shock. Furthermore, even though we demonstrate our main results for the nonlinear regularisation of the form~\cref{equation:nonlinear_fourth_order_regularisation_2_2}, one can also obtain analogous modified equal area rules using other regularised equations. For example, the \(\varepsilon \to 0\) limit of the equation
                \begin{equation}
                    \label{equation:nonlinear_fourth_order_regularisation_3_1}%
                    \pd{u}{t} = \Phi(u)_{xx} -\left(\varepsilon^2 f(u) \frac{\partial^3 u}{\partial x^3}\right)_{x} + R(u),       
                \end{equation}
                leads to shock position that satisfies an equal area rule in \(\int D(u)/f(u) \df{u}\), and this regularisation can lead to shock positions that do not conserve \(\Phi(u)\) across the shock. The nonlinear regularisation
                \begin{equation}
                    \label{equation:nonlinear_fourth_order_regularisation_0_4}%
                    \pd{u}{t} = \pdn{\Phi(u)}{x}{2} - \varepsilon^2 f(u)_{xxxx} + R(u),
                \end{equation}
                which is equivalent to
                \begin{equation}
                    \label{equation:nonlinear_fourth_order_regularisation_1_3}%
                    \pd{u}{t} = \pdn{\Phi(u)}{x}{2} - \varepsilon^2 \left(f'(u) \frac{\partial u}{\partial x}\right)_{xxx} + R(u),
                \end{equation}
                also leads to a different modified equal area rule
                \begin{equation}
                    \label{equation:modified_equal_area_rule_second_version}%
                    \int_{\ShockMin}^{\ShockMax} f'(u)(\Phi(u) - \Phi_S) \df{u} = 0,
                \end{equation}
                noting that for \Cref{equation:nonlinear_fourth_order_regularisation_0_4} both \(f(u)\) and \(f'(u)\) will need to be positive.

        \subsection{Numerical solutions obey theoretical shock dynamics}\label{section:numerical}
\begin{figure}[htbp!]
                \centering
                \subcaptionbox{\(D(u)\) \label{figure:numericalDiffusion}}{\includegraphics[width=0.49\linewidth]{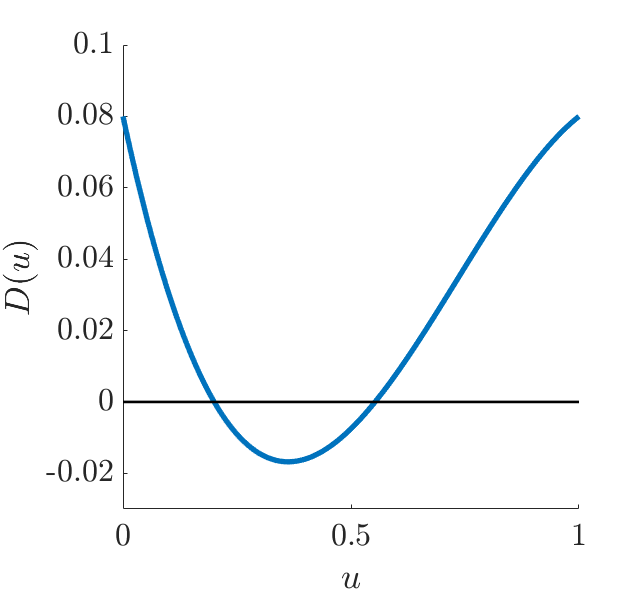}} 
                \subcaptionbox{\(R(u)\) \label{figure:numericalReaction}}{\includegraphics[width=0.49\linewidth]{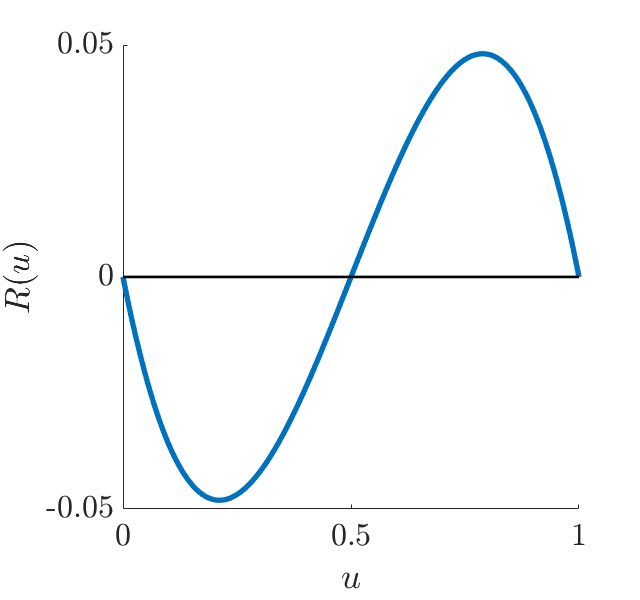}}
                \caption{Diffusivity and reaction term used for numerical solutions. (a) Diffusivity of the cubic form~\Cref{equation:quadratic_cubic_diffusion} with \(a = 0.2\), \(b = 0.4\) and \(\delta = 0.5\). (b) Reaction term of the cubic form~\Cref{equation:cubic_reaction} with \(\gamma = 0.5\).}
                \label{figures:numericalDiffusion+numericalReaction}
            \end{figure}

            We have found a nonlinear regularisation \Cref{equation:nonlinear_fourth_order_regularisation_2_2} which we can use to obtain a continuous diffusivity shock. To demonstrate numerically how a regularisation selects a shock position, we solve both the reaction--diffusion equation with the linear regularisation~\Cref{equation:standard_fourth_order_regularisation} and the nonlinear regularisation~\cref{equation:nonlinear_fourth_order_regularisation_2_2} using the method of lines. In both cases we used the cubic diffusivity~\Cref{equation:quadratic_cubic_diffusion} with \(a = 0.2\), \(b = 0.4\), and \(\delta = 0.5\) (see~\Cref{figure:numericalDiffusion}), which gives a relatively large spacing between the equal area in \(\Phi(u)\) shock position and the continuous diffusivity shock position. For the reaction term we used the cubic reaction~\Cref{equation:cubic_reaction} with \(\gamma = 0.5\), giving \(R(u) = u(1 - u)(u-0.5)\), as shown in \Cref{figure:numericalReaction}. When using the nonlinear regularisation~\cref{equation:nonlinear_fourth_order_regularisation_2_2}, we chose \(f(u) = \e^{-Au}\) with \(A = -3.0757\) (see~\Cref{figure:exponential_f}), which we previously showed gives the continuous diffusivity shock. We discretise space using central differences, and use a Heaviside step function as the initial condition. The uniform grid spacing was \(\Delta x = 0.001\), and we solved on a grid from \(x = 0\) to \(x = 10\) with Dirichlet boundary conditions \(x(0)= 1,\) \(x(10) = 0\), and \(x'(0)=x'(10)=0\), from \(t = 0\) to \(t = 20\).
            
\begin{figure}[htbp!]
                \centering
                \subcaptionbox{\label{figure:profiles_equal_area_not_zoomed}}{\includegraphics[width=0.49\linewidth]{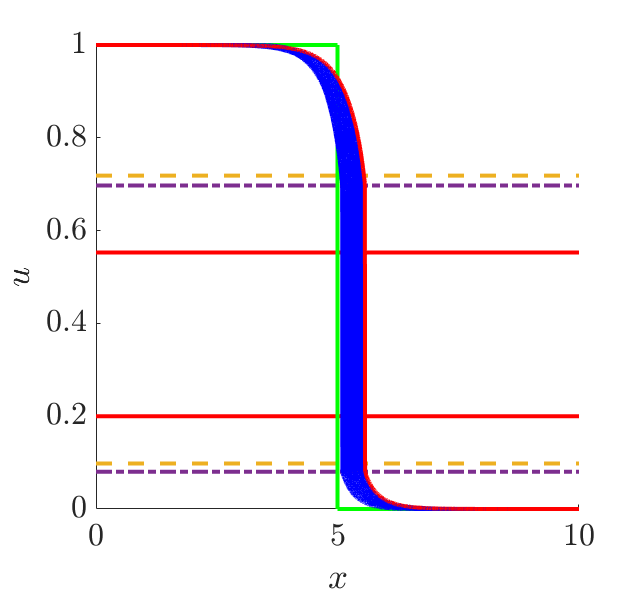}} 
                \subcaptionbox{\label{figure:profiles_equal_area_zoomed}}{\includegraphics[width=0.49\linewidth]{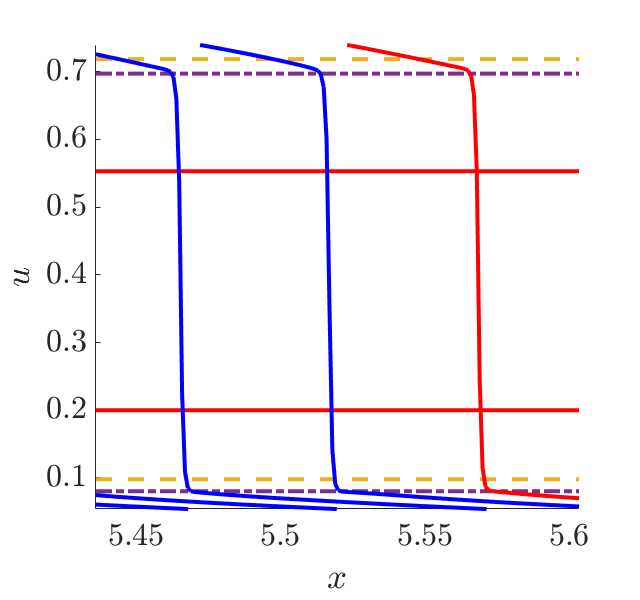}}   
                \subcaptionbox{\label{figure:profiles_continous_diffusion_not_zoomed}}{\includegraphics[width=0.49\linewidth]{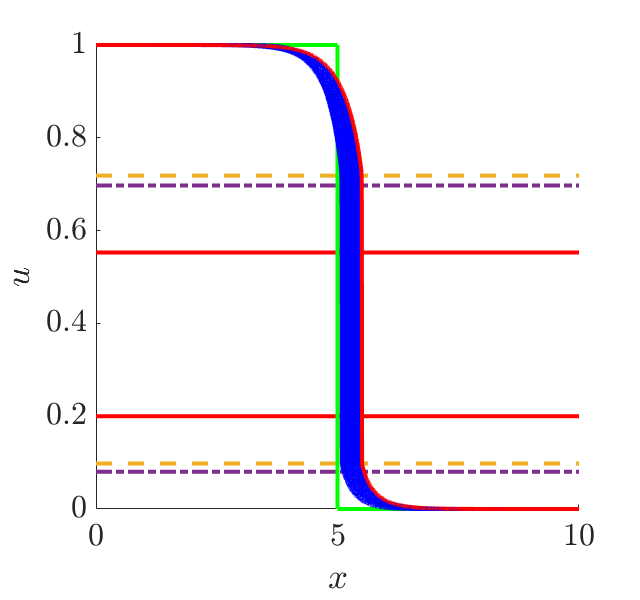}} 
                \subcaptionbox{\label{figure:profiles_continous_diffusion_zoomed}}{\includegraphics[width=0.49\linewidth]{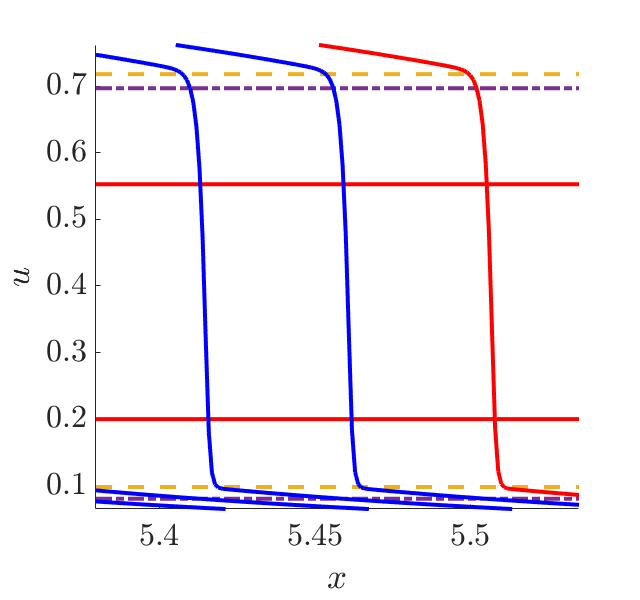}} 
                \caption{Numerical solutions with the diffusivity and reaction term shown in~\Cref{figures:numericalDiffusion+numericalReaction}. Green curves: Heaviside function initial condition. Red curve: Final profiles. Blue curves: Solutions at intermediate times. Solid red horizontal lines: zeros of the diffusivity. Dashed yellow lines: end-points of the continuous diffusivity shock. Dashed purple lines: end-points of the equal area in \(\Phi(u)\) shock. (a) Numerical solution of the reaction--diffusion equation with linear non-local regularisation~\Cref{equation:standard_fourth_order_regularisation} (b) Zoomed view of the last three profiles of (a). The shock appears to correspond to the equal area in \(\Phi(u)\) shock end-points (purple curve). (c) Numerical solution of the reaction--diffusion equation with nonlinear regularisation~\Cref{equation:nonlinear_fourth_order_regularisation_2_2}, where \(f(u) = \e^{-Au}\) with \(A = -3.0757\). (d) Zoomed view of the last three profiles of (c). The shock appears to correspond to the continuous diffusivity shock end-points (yellow curve).}
                \label{figures:numerical_profiles}
            \end{figure}

            \Cref{figure:profiles_equal_area_not_zoomed,figure:profiles_equal_area_zoomed} show solutions obtained with the linear regularisation~\Cref{equation:standard_fourth_order_regularisation} that~\textcite{witelski1995shocks,li2021shock} showed yields the equal area in \(\Phi(u)\) rule. The steepest part of the solution, or the shock, starts and ends approximately at the equal area in \(\Phi(u)\) shock end-points. \Cref{figure:profiles_continous_diffusion_not_zoomed,figure:profiles_continous_diffusion_zoomed} shows the profiles after solving the reaction--diffusion equation with nonlinear regularisation~\Cref{equation:nonlinear_fourth_order_regularisation_2_2}. As predicted, the shock now starts and ends approximately at the continuous diffusivity shock end-points.
    
\begin{figure}[htbp!]
                \centering
                \subcaptionbox{\label{figure:profiles_comparison}}{\includegraphics[width=0.45\linewidth]{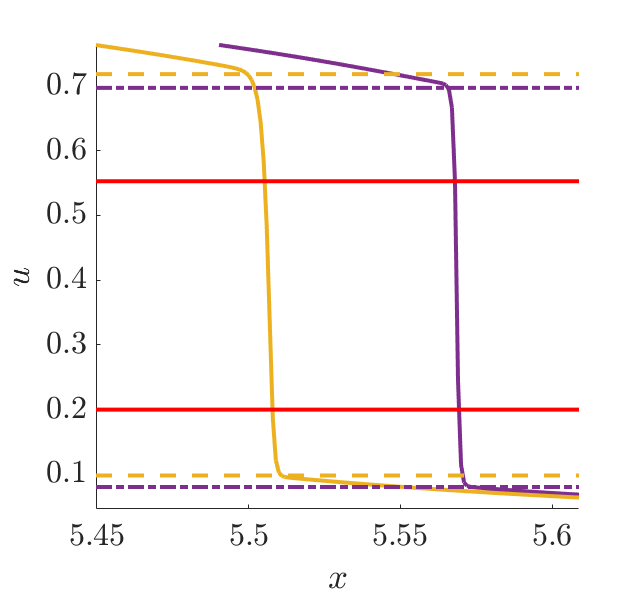}}
                \subcaptionbox{\label{figure:speed_comparison}}{\includegraphics[width=0.45\linewidth]{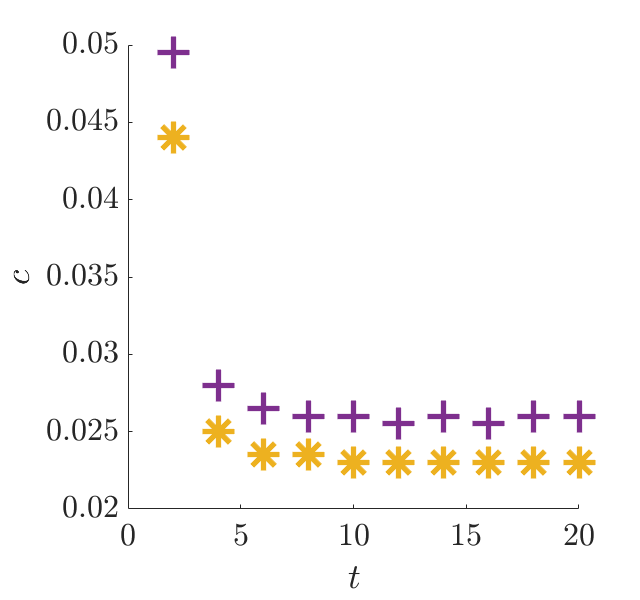}}
                \subcaptionbox{\label{figure:phase_plane_equal_area}}{\includegraphics[width=0.45\linewidth]{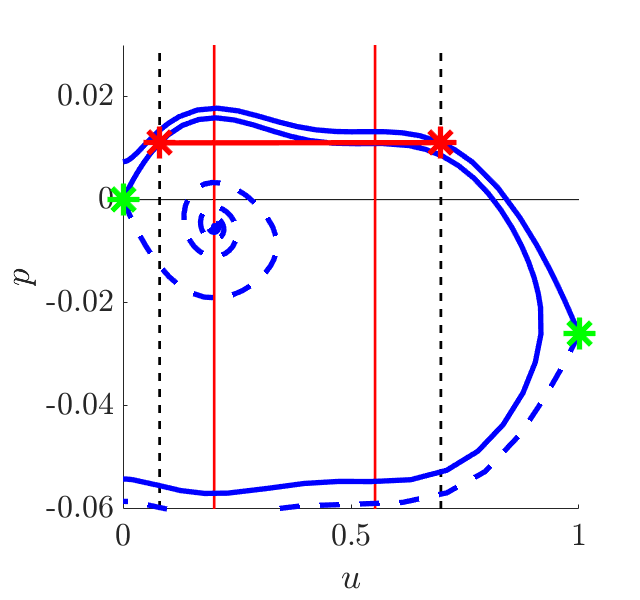}} 
                \subcaptionbox{\label{figure:phase_plane_continous_diffusion}}{\includegraphics[width=0.45\linewidth]{phase_plane_equal_diffusion_shock.png}}    
                \caption{Comparison of the wave speed in equal area in \(\Phi(u)\) and continuous diffusivity shocks. (a) Final profile of \Cref{figure:profiles_equal_area_not_zoomed,figure:profiles_continous_diffusion_not_zoomed} as the purple (equal area in \(\Phi(u)\)) and yellow (continuous diffusivity) curves, respectively. (b) The estimated speeds form the numerical solutions shown in~\Cref{figures:numerical_profiles}. Purple stars: estimated speeds for the solutions shown in \Cref{figure:profiles_equal_area_not_zoomed}. Yellow crosses: estimated speeds for the solutions shown in \Cref{figure:profiles_continous_diffusion_not_zoomed}. (c)--(d) Phase planes of the reduced problem~\Cref{equations:phase_plane}, similar to~\Cref{figures:phase_planes}. (c) The vertical dashed black lines are fixed according to the equal area in \(\Phi(u)\) shock position, with corresponding wave speed \(c = 0.026\). (d) The vertical dashed black lines are fixed according to the continuous diffusivity shock, with corresponding wave speed \(c = 0.0232\).}       
                \label{figures:phase_planes_numerical}
            \end{figure}

            \Cref{figure:profiles_comparison} shows the final profile from both numerical solutions. In addition to moving the shock position, the different regularisations also influence the wave speed. Using the numerical solutions, we estimate the speed of the solution by locating the \(x\)-position of the shock at each time step. We define the shock \(x\)-position to be the position \(x\) where there is the largest value of the slope \(|\uppartial u/\uppartial x|\) estimated using the central difference. To obtain the wave speed, we find the shock \(x\)-position at \(t \in \{0,2,4,6,8,10,12,14,16,18,20\},\) and obtain backward difference estimates for the speed at \(t \in \{2,4,6,8,10,12,14,16,18,20\}.\) \Cref{figure:speed_comparison} compares the estimated speeds of the two solutions. Since the solutions evolve towards travelling waves, the solutions quickly settle to roughly constant speeds. The estimated speed in the equal area in \(\Phi(u)\) shock solution is \(c \approx 0.026,\) whereas and the estimated speed of continuous diffusivity shock solution is \(c \approx 0.023\).  These differences in speed can be explained using the phase plane of the reduced problem~\Cref{equations:phase_plane}, which is the same for both regularised equations. \Cref{figures:phase_planes_numerical} show the phase plane for both shocks. A shock is only possible if the value of \(p\) is constant at \(u = \ShockMin\) and \(u = \ShockMax\). At our final estimated speeds this constant condition is satisfied, confirming that our numerical solutions are following the theoretical shock dynamics we derived previously.
    
    \section{Conclusion}
        In this paper we used a nonlinear regularisation in order to control the shock position and obtain a continuous diffusivity shock in a reaction--diffusion equation with partially negative diffusivity. This regularisation gave a modified equal area rule \Cref{equation:modified_equal_area_rule} that by choice of regularisation function \(f(u)\) could control the shock position. The continuous diffusivity shock position is only distinct from the equal area in \(\Phi(u)\) shock in the case of non-symmetric diffusivity. We showed that the continuous diffusivity rule for decreasing-increasing diffusivities selects the longest possible shock.  We explored these dynamics and used geometric singular perturbation theory to prove the existence of travelling wave solutions with the continuous diffusivity shock for cubic reaction terms. The nonlinear regularisation and modified equal area rule derived in this paper also apply to arbitrary reaction terms, including the nonlinear diffusion equation where \(R(u) = 0\).

        In our numerical solutions we solved reaction--diffusion equations with linear non-local~\Cref{equation:standard_fourth_order_regularisation} and nonlinear~\cref{equation:nonlinear_fourth_order_regularisation_2_2} regularisations. In the absence of a regularisation imposed explicitly, numerical error in the discretisation can also implicitly regularise the reaction--diffusion equation. For our numerical solutions, we used the central difference
        \begin{equation}
            \pdn{\Phi(u)}{x}{2} \approx \frac{\Phi_{i - 1} - 2\Phi_i + \Phi_{i + 1}}{(\Delta x)^2}, 
            \label{equation:our_numerical_discretisation}
        \end{equation}
        where \(\Delta x\) is the grid spacing, which to leading order has the error
        \begin{equation}
            \label{equation:our_numerical_error}%
            E = -(\Delta x)^2\left[\frac{\Phi_{xxxx}}{12}\right]_{i} + \mathcal{O}(\Delta x^4)
        \end{equation}
        as \(\Delta x \to 0,\) which yields a regularisation of the form~\Cref{equation:nonlinear_fourth_order_regularisation_0_4}, for which the corresponding shock condition is~\Cref{equation:modified_equal_area_rule_second_version}. Using this rule to anslyse our numerical error term we get
        \begin{equation}
            \label{equation:our_numerical_error_shock_rule}%
            \int_{\ShockMin}^{\ShockMax} D(u)\left(\Phi(u) - \Phi_S\right)\df{u} = \left[\frac{\Phi(u)^2}{2} - \Phi(u)\Phi_S\right]_{u = \ShockMin}^{u = \ShockMax},  
        \end{equation}
        which is satisfied by any shock with continuous \(\Phi(u)\). Alternative discretisations of the diffusivity are possible, for example the conservation form
        \begin{equation}
            \label{equation:other_numerical_discretisation}%
            \pdn{\Phi(u)}{x}{2} = \pd{}{x}\left(D(u)\pd{u}{x}\right) \approx \frac{D_{i+1/2}(u_{i + 1} - u_i) - D_{i - 1/2}(u_i - u_{i - 1})}{(\Delta x)^2},   
        \end{equation}
        where
        \begin{equation}
            \label{equation:other_numerical_discretisation_part_two}%
            D_{i+1/2} = \frac{D_i + D_{i + 1}}{2}.
        \end{equation}
        This discretisation has the numerical error
        \begin{equation}
            \label{equation:other_numerical_error}%
            E = -\frac{(\Delta x)^2}{12}\left[Du_{xxxx} + 2D_xu_{xxx} + 3D_{xx}u_{xx} + 2D_{xxx}u_{x}\right]_i + \mathcal{O}(\Delta x^4)
        \end{equation}
        as \(\Delta x \to 0.\) This form of the numerical error corresponds to the regularisation
        \begin{equation}
            \label{equation:other_numerical_error_regularisation}%
            -\varepsilon^2 \left[\pdn{\Phi(u)}{x}{4} + \pd{}{x}\left(D''(u)\left(\pd{u}{x}\right)^3\right)\right],    
        \end{equation}
        for which we do not have a shock selection rule. Preliminary numerical solutions of the regularised reaction--diffusion equations \Cref{equation:standard_fourth_order_regularisation,equation:nonlinear_fourth_order_regularisation_2_2} using the discretisation \Cref{equation:other_numerical_discretisation} seem to show the shocks being moved away from their predicted shock positions. A possible avenue of future research would be to find the shock condition corresponding to the regularisation term \Cref{equation:other_numerical_error_regularisation} and to see if using the discretisation \Cref{equation:other_numerical_discretisation} is acting to move the shocks to ones that satisfy that shock condition.      
        
        Another potential extension is the consider different forms of the diffusivity \(D(u).\) In this paper we focus on positive-negative-positive diffusivities where the diffusivity has two zeros, but this is not the only possibility. For instance~\textcite{maini2006diffusion,maini2007aggregative} considered reaction--diffusion equations in the case where the diffusivity had one zero and one sign change. \textcite{kuzmin2011front} mainly considered a positive-negative-diffusivity but generalised their results to diffusivities with \(2n\) sign changes. Given that the continuous diffusivity shock is of the longest possible length, another avenue for future work is to identify the consequences for this solution in terms of the entropy of the solutions. Given that there are many other possible shocks, entropy considerations might also provide another way to select a shock position.
    \section*{Declaration of competing interests}
        The authors declare that they have no known competing financial interests or personal relationships that could have appeared to influence the work reported in this paper.
    \section*{Acknowledgements}
        This work was supported by funding from the Australian Research Council Discovery Program (grant numbers DP200102130, DP230100406, DE240100897).
    \printbibliography    
\end{document}